\documentclass[12pt]{article}
\usepackage{amsmath}    
\usepackage{amssymb}    
\usepackage{latexsym}
\usepackage{graphicx}
\usepackage{multirow}
\newcommand{\Pn}    {\mathbb{P}}
\newcommand{\PP}    {{\bf P}}
\newcommand{\e}     {\ensuremath{\varepsilon}}
\newcommand{\dr}    {\ensuremath{\partial }}

\newcommand{\transp}[1] {#1 ^{\top} }

\newcommand{\R} {\mathbb{R}}

\newcommand{\E} {\mathbb{E}}

\newcommand{\cA}{{\mathcal A}}
\newcommand{\cB}{{\mathcal B}}
\newcommand{\cC}{{\mathcal C}}
\newcommand{\cD}{{\mathcal D}}
\newcommand{\cE}{{\mathcal E}}
\newcommand{\cF}{{\mathcal F}}
\newcommand{\cG}{{\mathcal G}}
\newcommand{\cH}{{\mathcal H}}

\newcommand{\cN}{{\mathcal N}}

\newcommand{\cR}{{\mathcal R}}
\newcommand{\cS}{{\mathcal S}}

\newcommand{\cvL }       { \stackrel   { \it \mathcal Law  } { \rightarrow} }

\newtheorem{definition}  {\bf Definition}
\newtheorem{theorem}     {\bf Theorem}
\newtheorem{lemma}       {\bf Lemma}
\newtheorem{corollary}   {\bf Corollary}
\newtheorem{proposition} {\bf Proposition}

\newtheorem{property}    {\bf Property}
\newtheorem{example}     {$\mathcal Example$}
\newtheorem{Sim}         {$\mathcal Simulation$}
\newtheorem{remark}      {$\mathcal Remark$ }

\begin{document}
\centerline{\noindent \bf \Large PROJECTION PURSUIT THROUGH}
\centerline{\noindent \bf \Large RELATIVE ENTROPY MINIMISATION}
\vskip 3mm

\vskip 5mm
{\large\center\noindent  Jacques Touboul

\noindent Laboratoire de Statistique Th\'eorique et Appliqu\'ee

\noindent Universit\'e Pierre et Marie Curie

\hskip 6cm \noindent jack\_touboul@hotmail.com
}

\vskip 5mm
{\it Projection Pursuit methodology permits to solve the difficult problem of finding an estimate of a density defined on a set of very large dimension.
In his seminal article, Huber (see "Projection pursuit", Annals of Statistics, 1985) evidences the interest of the Projection Pursuit method thanks to the factorisation of a density into a Gaussian component and some residual density in a context of Kullback-Leibler divergence maximisation.\\
In the present article, we introduce a new algorithm, and in particular  a test for the factorisation of a density estimated from an iid sample.\\}
\vskip 4mm
\noindent {\bf1. Outline of the article}

Projection Pursuit aims at creating one or several projections delivering a maximum of information on the structure of a data set irrespective of its size. Once a structure has been evidenced, the corresponding data are transformed through a Gaussianization. Recursively, this process is repeated in order to determine another structure in the remaining data until no further structure can be highlighted eventually. These kind of approaches for isolating structures were first studied by Friedman \cite{Frie84} and Huber \cite{MR790553}. Each of them details, through two different methodologies each, how to isolate such a structure and therefore how to estimate the density of the corresponding data.\\
However, since Mu Zhu \cite{Zhu2004} showed the two methodologies described by each of the above authors did not in fact turn out to be equivalent when the number of iterations in the algorithms exceeds the dimension of the space containing the data, we will consequently only concentrate on Huber's study while taking into account Mu Zhu's input.

After providing a brief overview of Huber's methodologies, we will then expose our approach and objective.
\vskip 3mm
\noindent {\bf1.1.}{\it Huber's analytic approach}

A density $f$ on $\R^d$ is considered. We then define an instrumental density $g$ with the same mean and variance as  $f$.
According to Huber's approach, we first carry out the $K(f,g)=0$ test - with $K$ being the relative entropy (also called the Kullback-Leibler divergence). If the test is passed, then $f=g$ and the algorithm stops. If the test were not to be verified, based on the maximisation of $a\mapsto K(f_a,g_a)$ since $K(f,g)=K(f_a,g_a)+K(f\frac{g_a}{f_a},g)$ and assuming that $K(f,g)$ is finite, Huber's methodology requires as a first step to define a vector $a_1$ and a density $f^{(1)}$  with
\begin{equation}\label{DefSequMethodH1}
a_1\ =\ arg\inf_{a\in\R^d_*}\ K(f\frac{g_a}{f_a},g)\text{ and }f^{(1)}=f\frac{g_{a_1}}{f_{a_1}},
\end{equation}
where $\R^d_*$ is the set of non null vectors of $\R^d$ and $f_a$ (resp. $g_a$) represents the density of $\transp aX$ (resp. $\transp aY$) when $f$ (resp. $g$) is the density of $X$ (resp. $Y$). \\
As a second step, Huber's algorithm replaces  $f$ with $f^{(1)}$ and repeats the first step.\\
Finally, a sequence $(a_1,a_2,\ldots)$ of vectors of $\R^d_*$ and a sequence of densities $f^{(i)}$ are derived from the iterations of this process.
\begin{remark}
$\\$The algorithm enables us to generate a product approximation and even a product representation of $f$. Indeed, two rules can trigger the end of the process. The first one is the nullity of the relative entropy and the second one is the process reaching the $d^{th}$ iteration. When these two rules are satisfied, the algorithm produces a product approximation of $f$. When only the first rule is satisfied, the algorithm generates a product representation of $f$.\\
Mathematically, for any integer $j$, such that $K(f^{(j)},g)=0$ with $j\leq d$, the process infers $f^{(j)}=g$, i.e. $f=g\Pi_{i=1}^j\frac{f^{(i-1)}_{a_i}} {g_{a_i}}$ since by induction $f^{(j)}=f\Pi_{i=1}^j\frac{g_{a_i}}{f^{(i-1)}_{a_i}}$. Likewise, when, for all $j$, it gets $K(f^{(j)},g)>0$ with $j\leq d$, it is assumed $g=f^{(d)}$ in order to obtain $f=g\Pi_{i=1}^d\frac{f^{(i-1)}_{a_i}} {g_{a_i}}$, i.e. we approximate $f$ with the product $g\Pi_{i=1}^d\frac{f^{(i-1)}_{a_i}} {g_{a_i}}$.\\
Even if the condition $j\leq d$ is not met, the algorithm can also stop if the Kullback-Leibler divergence equals zero. Therefore, since by induction we have $f^{(j)}=f\Pi_{i=1}^j\frac{g_{a_i}}{f^{(i-1)}_{a_i}}$ with $f^{(0)}=f$, we infer $g=f\Pi_{i=1}^j\frac{g_{a_i}}{f^{(i-1)}_{a_i}}$. We can thus represent $f$ as $f=g\Pi_{i=1}^j\frac{f^{(i-1)}_{a_i}} {g_{a_i}}.$\\
Finally, we remark that the algorithm implies that the sequence $(K(f^{(j)},g))_j$ is decreasing and non negative with $f^{(0)}=f$.
\end{remark}
\vskip 3mm
\noindent {\bf1.2.}{\it  Huber's synthetic approach}

Maintaining the notations of the above section, we begin with performing the $K(f,g)=0$ test; If the test is passed, then $f=g$ and the algorithm stops, otherwise, based on the maximisation of $a\mapsto K(f_a,g_a)$ since $K(f,g)=K(f_a,g_a)+K(f,g\frac{f_a}{g_a})$ and assuming that $K(f,g)$ is finite, Huber's methodology requires as a first step to define a vector $a_1$ and a density $g^{(1)}$  with
\begin{equation}\label{DefSequMethodH2}
a_1\ =\ arg\inf_{a\in\R^d_*}\ K(f,g\frac{f_a}{g_a})\text{ and }g^{(1)}=g\frac{f_{a_1}}{g_{a_1}}.
\end{equation}
As a second step, Huber's algorithm replaces  $g$ with $g^{(1)}$ and repeats the first step.\\
Finally, a sequence $(a_1,a_2,\ldots)$ of vectors of $\R^d_*$ and a sequence of densities $g^{(i)}$  are derived from the iterations of this process.
\begin{remark}
$\\$Similarly as in the analytic approach, this methodology allows us to generate a product approximation and even a product representation of $f$ from $g$. Moreover, it also offers the same end of process rules.\\
In other words, if for any $j$, such that $j\leq d$, we have $K(f,g^{(j)})>0$, then $f$ is approximated with $g^{(d)}$. And if there exists $j$, such that  $K(f,g^{(j)})=0$, then it holds $g^{(j)}=f$, i.e. $f$ is represented by $g^{(j)}$. In this case, the relationship $K(f,g^{(j)})=0$ implies that  $g^{(j)}=f$, i.e. since by induction we have $g^{(j)}=g\Pi_{i=1}^j\frac{f_{a_i}}{g^{(i-1)}_{a_i}}$ with $g^{(0)}=g$, it holds $f=g\Pi_{i=1}^j\frac{f_{a_i}}{g^{(i-1)}_{a_i}}$.\\
Eventually, we note that the algorithm implies that the sequence $(K(f,g^{(j)}))_j$ is decreasing and non negative with  $g^{(0)}=g$.
\end{remark}
\noindent Finally, in \cite{Zhu2004}, Mu Zhu shows that, beyond $d$ iterations, the data processing of these methodologies evidences significant differences, i.e. that past $d$ iterations, the two methodologies are no longer equivalent. We will therefore only consider Huber's synthetic approach since  $g$ is known and since we want to find a representation of $f$.
\vskip 3mm
\noindent {\bf1.3.}{\it Proposal}\label{OurProposal}

We begin with performing the $K(f,g)=0$ test; should this test be passed, then $f=g$ and the algorithm stops, otherwise, the first step of our algorithm consists in defining a vector $a_1$ and a density $g^{(1)}$  by 
\begin{equation}\label{DefSequMethod}
a_1\ =\ arg\inf_{a\in\R^d_*}\ K(g\frac{f_a}{g_a},f)\text{ and }g^{(1)}=g\frac{f_{a_1}}{g_{a_1}}.
\end{equation}
In the second step, we replace $g$ with $g^{(1)}$ and we repeat the first step. We thus derive, from the iterations of this process, a sequence  $(a_1,a_2,...)$ of vectors in  $\R^d_*$  and a sequence of densities $g^{(i)}$.
We will prove that $a_1$ simultaneously optimises (\ref{DefSequMethodH1}), (\ref{DefSequMethodH2})  and (\ref{DefSequMethod}). We will also prove that the underlying structures of $f$ evidenced through this method are identical to the ones obtained through the Huber's methods.
\begin{remark}
$\\$As in Huber's algorithms, we perform a product approximation and even a product representation of $f$.\\
In the case where, at each of the $d^{th}$ first steps, the relative entropy is positive, we then  approximate $f$ with $g^{(d)}$.\\
In the case where there exists a step of the algorithm such that the Kullback-Leibler divergence equals zero, then, calling $j$ this step, we represent $f$ with $g^{(dj}$. In other words, if there exists a positive integer $j$ such that $K(g^{(j)},f)=0$, then, since by induction we have $g^{(j)}=g\Pi_{i=1}^j\frac{f_{a_i}}{g^{(i-1)}_{a_i}}$ with $g^{(0)}=g$, we represent $f$ with the product $g\Pi_{i=1}^j\frac{f_{a_i}}{g^{(i-1)}_{a_i}}.$\\
We also remark that the algorithm implies that the sequence $(K(g^{(j)},f))_j$ is decreasing and non negative with  $g^{(0)}=g$.\\
Finally, the very form of the relationship (\ref{DefSequMethod}) demonstrates that we deal with M-estimation. We can consequently state that our method is more robust than Huber's - see \cite{BAR-YOHAI}, \cite{TOMA} as well as \cite{RobStat}.
\end{remark}
\begin{example}Let $f$ be a density defined on $\R^{10}$ by $f(x_1,\ldots,x_{10})=\eta(x_2,\ldots,x_{10})\zeta(x_1)$, with $\eta$ being a multivariate Gaussian density on $\R^{9}$, and $\zeta$ being a non Gaussian density.\\
Let us also consider $g$, a multivariate Gaussian density with the same mean and variance as $f$.\\
Since $g(x_2,\ldots,x_{10}/x_1)=\eta(x_2,\ldots,x_{10})$, we have 
$K(g\frac{f_1}{g_1},f)=K(\eta.f_1,f)=K(f,f)=0$ as $f_1=\zeta$ - where $f_1$ and $g_1$ are the first marginal
densities of $f$ and $g$ respectively. Hence, the non negative function $a\mapsto K(g\frac{f_a}{g_a},f)$ reaches zero for $e_1=(1,0,\ldots,0)'$. \\
We therefore infer that $g(x_2,\ldots,x_{10}/x_1)=f(x_2,\ldots,x_{10}/x_1).$
\end{example}
To recapitulate our method, if $K(g,f)=0$, we derive $f$ from the relationship $f=g$; should a sequence $(a_i)_{i=1,...j}$, $j<d$, of vectors in $\R^d_*$ defining $g^{(j)}$ and such that $K(g^{(j)},f)=0$ exist, then $f(./\transp{a_i}x, 1\leq i\leq j)=g(./\transp{a_i}x, 1\leq i\leq j)$, i.e. $f$ coincides with $g$ on the complement  of the vector subspace generated by the family $\{a_i\}_{i=1,...,j}$ - see also section 2.1.2. for details.

In this paper, after having clarified the choice of $g$,  we will consider the statistical solution to the representation problem, assuming that $f$ is unknown and $X_1$, $X_2$,... $X_m$ are i.i.d. with density $f$. We will provide asymptotic results pertaining to the family of optimizing vectors $a_{k,m}$ - that we will define more precisely below - as $m$ goes to infinity.
Our results also prove that the empirical representation scheme converges towards the theoretical one.
Finally, we will compare  Huber's optimisation methods with ours and we will present simulations.
\vskip 3mm
\noindent {\bf2.  The algorithm}\label{TheAlgo}\\
\noindent {\bf2.1.}{\it  The model}\label{modelSection}

As described by Friedman \cite{Frie84} and Diaconis \cite{MR0751274}, the choice of $g$ depends on the family of distribution one wants to find in $f$. Until now, the choice has only been to use the class of Gaussian distributions. This can also be extended to the class of elliptical distributions.
\vskip 3mm
\noindent {\bf2.1.1.}{\it  Elliptical distributions}

The fact that conditional densities with elliptical distributions are also elliptical - see \cite{MR0629795}, \cite{MR2061237} - enables us to use this class in our algorithm - and in Huber's algorithms.
\begin{definition}
$X$ is said to abide by a multivariate elliptical distribution, denoted $X\sim E_d(\mu,\Sigma,\xi_d)$, if $X$ has the following density, for any $x$ in $\R^d$ : 
$f_X(x)=\frac{c_d}{|\Sigma|^{1/2}}\xi_d\Big(\frac{1}{2}(x-\mu)'\Sigma^{-1}(x-\mu)\Big)$,\\
where $\Sigma$ is a $d\times d$ positive-definite matrix and where $\mu$ is an $ d$-column vector,\\
where $\xi_d$ is referred as the "density generator",\\
where $c_d$ is a normalisation constant, such that
$c_d=\frac{\Gamma(d/2)}{(2\pi)^{d/2}}\Big(\int_0^\infty x^{d/2-1}\xi_d(x)dx\Big)^{-1}$,\\ with $\int_0^\infty x^{d/2-1}\xi_d(x)dx<\infty$.
\end{definition}
\begin{property}\label{ElliptProp}
1/ For any $X\sim E_d(\mu,\Sigma,\xi_d)$, for any $m\times d$ matrix with rank $m\leq d,$ $A$, and for any  $m$-dimensional vector, $b$, we have $AX+b\sim E_m(A\mu+b,A\Sigma A',\xi_m)$.\\
Any marginal density of multivarite elliptical distribution is consequently elliptical, i.e. \\
$X=(X_1,X_2,...,X_d)\sim E_d(\mu,\Sigma,\xi_d)$ implies that $X_i\sim E_1(\mu_i,\sigma^2_i,\xi_1)$ with  $f_{X_i}(x)= \frac{c_1}{\sigma_i}\xi_1\Big(\frac{1}{2}(\frac{x-\mu_i}{\sigma})^2\Big),$ $1\leq i\leq d$.\\
2/ Corollary 5 of \cite{MR0629795} states that the conditional densities with elliptical distributions are also elliptical. Indeed, if $X=(X_1,X_2)'\sim E_d(\mu,\Sigma,\xi_d)$, with $X_1$ (resp. $X_2$) of size $d_1<d$ (resp. $d_2<d$), then $X_1/(X_2=a)\sim E_{d_1}(\mu',\Sigma',\xi_{d_1})$ with $\mu'=\mu_1+\Sigma_{12}\Sigma_{22}^{-1}(a-\mu_2)$ and $\Sigma'=\Sigma_{11}-\Sigma_{12}\Sigma_{22}^{-1}\Sigma_{21},$
with $\mu=(\mu_1,\mu_2)$ and $\Sigma=(\Sigma_{ij})_{1\leq i,j\leq 2}$.
\end{property}
\begin{remark}\label{implyEstimBounded}
In \cite{MR2061237}, the authors show that the multivariate Gaussian distribution derives from $\xi_d(x)=e^{-x}$. They also show that if $X=(X_1,...,X_d)$ has an elliptical density such that its marginals meet $E(X_i)<\infty$ and $E(X_i^2)<\infty$ for $1\leq i\leq d,$ then $\mu$ is the mean of $X$ and $\Sigma$ is a multiple of the covariance matrix of $X$. From now on, we will therefore assume this is the case.
\end{remark}
\begin{definition}
Let $t$ be an elliptical density on $\R^k$ and let $q$ be an elliptical density on $\R^{k'}$.
The elliptical densities $t$ and $q$ are said to be part of the same family of elliptical densities, if their generating densities are $\xi_k$ and $\xi_{k'}$ respectively, which belong to a common given family of densities.
\end{definition}
\begin{example}
Consider two Gaussian densities  $\cN(0,1)$ and $\cN((0,0),Id_2)$. They are said to belong to the same elliptical family as they both present  $x\mapsto e^{-x}$ as generating density.
\end{example}
\noindent {\bf2.1.2.}{\it Choice of $g$}\label{gChoice}

Let $f$ be a density on $\R^d$. We assume there exists $d$ non null linearly independent vectors $a_j$, with $1\leq j\leq d,$ of $\R^d$, such that 
\begin{equation}\label{f-def}
f(x)=n(\transp{a_{j+1}}x,...,\transp{a_{d}}x)h(\transp{a_{1}}x,...,\transp{a_{j}}x),
\end{equation}
with $j<d$, $n$ being an elliptical density on $\R^{d-j-1}$ and with $h$ being a density on $\R^{j}$, which does not belong to the same family as $n$. Let $X=(X_{1},...,X_{d})$ be a vector with $f$ as density.\\
We define $g$ as an elliptical distribution with the same mean and variance as $f$.\\
For simplicity, let us assume that the family $\{a_j\}_{1\leq j\leq d}$ is the canonical basis of $\R^d$:\\
The very definition of $f$ implies that $(X_{j+1},...,X_{d})$ is independent from $(X_{1},...,X_{j})$. Hence, the property \ref{ElliptProp} allows us to infer that the density of $(X_{j+1},...,X_{d})$ given $(X_{1},...,X_{j})$ is $n$.\\
Let us assume that $K(g^{(j)},f)=0,$ for some $j\leq d$. We then get $\frac{f(x)}{f_{a_1}f_{a_2}...f_{a_j}}=\frac{g(x)}{g^{(1-1)}_{a_1}g^{(2-1)}_{a_2}...g^{(j-1)}_{a_j}}$, since,  by induction, we have $g^{(j)}(x)=g(x)\frac{f_{a_1}}{g^{(1-1)}_{a_1}}\frac{f_{a_2}}{g^{(2-1)}_{a_2}}...\frac{f_{a_j}}{g^{(j-1)}_{a_j}}$.
Consequently, the fact that the conditional densities with elliptical distributions are also elliptical, as well as the above relationship enable us to state that
$n(\transp{a_{j+1}}x,.,\transp{a_{d}}x)=f(./\transp{a_i}x, 1\leq i\leq j)=g(./\transp{a_i}x, 1\leq i\leq j).$ In other words, $f$ coincides with $g$ on the complement  of the vector subspace generated by the family $\{a_i\}_{i=1,...,j}$.

At present, if the family $\{a_j\}_{1\leq j\leq d}$ is no longer the canonical basis of $\R^d$, then this family is again a basis of $\R^d$. Hence, lemma \ref{ChangBasis} implies that
\begin{equation}\label{RelElli36}
g(./\transp{a_{1}}x,...,\transp{a_{j}}x)=n(\transp{a_{j+1}}x,...,\transp{a_{d}}x)=f(./\transp{a_{1}}x,...,\transp{a_{j}}x),
\end{equation}
which is equivalent to $K(g^{(j)},f)=0$, since by induction  $g^{(j)}=g\frac{f_{a_1}}{g^{(1-1)}_{a_1}}\frac{f_{a_2}}{g^{(2-1)}_{a_2}}...\frac{f_{a_j}}{g^{(j-1)}_{a_j}}$.\\
The end of our algorithm implies that $f$ coincides with $g$ on the complement  of the vector subspace generated by the family $\{a_i\}_{i=1,...,j}$.
Therefore, the nullity of the Kullback-Leibler divergence provides us with information on the density structure. In summary, the following proposition clarifies the choice of $g$ which depends on the family of distribution one wants to find in $f$ :
\begin{proposition}\label{Pb1}With the above notations, $K(g^{(j)},f)=0$ is equivalent to

$g(./\transp{a_{1}}x,...,\transp{a_{j}}x)=f(./\transp{a_{1}}x,...,\transp{a_{j}}x).$
\end{proposition}
More generally, the above proposition leads us to defining the co-support of $f$ as the vector space generated by the vectors $a_{1},...,a_{j}$.
\begin{definition}
Let $f$ be a density on $\R^d$. We define the co-vectors of $f$ as the sequence of vectors $a_{1},...,a_{j}$ which solves the problem $K(g^{(j)},f)=0$ where $g$ is an elliptical distribution with the same mean and variance as $f$.
We define the co-support of $f$ as the vector space generated by the vectors $a_{1},...,a_{j}$.
\end{definition}
\noindent {\bf2.2.}{\it Stochastic outline of the algorithm}\label{UseSample}

Let $X_1$, $X_2$,..,$X_m$ (resp. $Y_1$, $Y_2$,..,$Y_m$) be a sequence of $m$ independent random vectors with the same density $f$ (resp. $g$). 
As customary in nonparametric Kullback-Leibler optimizations, all estimates of $f$ and $f_a$, as well as all uses of Monte Carlo methods are being performed using subsamples $X_1$, $X_2$,..,$X_n$ and $Y_1$, $Y_2$,..,$Y_n$, extracted respectively from $X_1$, $X_2$,..,$X_m$ and $Y_1$, $Y_2$,..,$Y_m$, since the estimates are bounded below by some positive deterministic sequence $\theta_m$ (see Appendix B).\\
Let $\Pn_n$ be the empirical measure based on the subsample $X_1$, $X_2$,.,$X_n$. Let $f_n$ (resp.  $f_{a,n}$ for any $a$ in $\R^d_*$) be the kernel estimate of $f$ (resp.  $f_a$), which is built from $X_1$, $X_2$,..,$X_n$ (resp. $\transp aX_1$, $\transp aX_2$,..,$\transp aX_n$).\\
As defined in section 1.3, we introduce the following sequences $(a_{k})_{k\geq 1}$ and $(g^{(k)})_{k\geq 1}$:\\
$\bullet$ $a_{k}$ is a non null vector of $\R^d$ such that $a_{k}=arg\min_{a\in\R^d_*}, K(g^{(k-1)}\frac{f_a}{g^{(k-1)}_a},f)$,\\
$\bullet$ $g^{(k)}$ is the density such that $g^{(k)}=g^{(k-1)}\frac{f_{a_{k}}}{g^{(k-1)}_{a_{k}}}$ with $g^{(0)}=g$.\\
The stochastic setting up of the algorithm uses $f_n$ and $g_n^{(0)}=g$ instead of $f$ and $g^{(0)}=g$, since $g$ is known. Thus, at the first step, we build the vector $\check a_1$ which minimizes the Kullback-Leibler divergence between $f_n$ and $g\frac{f_{a,n}}{g_{a}}$ and which estimates $a_1$.\\
Proposition \ref{QuotientDonneLoi} and lemma \ref{toattaint} enable us to minimize the Kullback-Leibler divergence between $f_n$ and $g\frac{f_{a,n}}{g_{a}}$. Defining $\check a_1$ as the argument of this minimization, proposition \ref{KernelpConv2} shows us that this vector tends to $a_1$.\\
Finally, we define the density $\check g^{(1)}_m$ as $\check g^{(1)}_m=g\frac{f_{\check a_1,m}}{g_{\check a_1}}$ which estimates $g^{(1)}$ through theorem \ref{KernelKRessultatPricipal}.\\
Now, from the second step and as defined in section 1.3, the density $g^{(k-1)}$ is unknown. Once again, we therefore have to truncate the samples.\\
All estimates of $f$ and $f_a$ (resp. $g^{(1)}$ and $g_a^{(1)}$) are being performed using a subsample $X_1$, $X_2$,..,$X_n$ (resp. $Y_1^{(1)}$, $Y_2^{(1)}$,..,$Y_n^{(1)}$) extracted from $X_1$, $X_2$,..,$X_m$ (resp. $Y_1^{(1)}$, $Y_2^{(1)}$,..,$Y_m^{(1)}$ - which is a sequence of $m$ independent random vectors with the same density $g^{(1)}$) such that the estimates are bounded below by some positive deterministic sequence $\theta_m$ (see Appendix B).\\
Let $\Pn_n$ be the empirical measure based on the subsample $X_1$, $X_2$,..,$X_n$. Let $f_n$ (resp. $g_n^{(1)}$, $f_{a,n}$, $g_{a,n}^{(1)}$ for any $a$ in $\R^d_*$) be the kernel estimate of $f$ (resp. $g^{(1)}$, $f_a$, $g_a^{(1)}$) which is built from $X_1$, $X_2$,..,$X_n$ (resp. $Y_1^{(1)}$, $Y_2^{(1)}$,..,$Y_n^{(1)}$).
The stochastic setting up of the algorithm uses $f_n$ and $g_n^{(1)}$ instead of $f$ and $g^{(1)}$.
Thus, we build the vector $\check a_2$ which minimizes the Kullback-Leibler divergence between $f_n$ and $g_n^{(1)}\frac{f_{a,n}}{g_{a,n}^{(1)}}$ - since $g^{(1)}$ and $g^{(1)}_a$ are unknown - and which estimates $a_2$.
Proposition \ref{QuotientDonneLoi} and lemma \ref{toattaint} enable us to minimize the Kullback-Leibler divergence between $f_n$ and $g_n^{(1)}\frac{f_{a,n}}{g_{a,n}^{(1)}}$. Defining $\check a_2$ as the argument of this minimization, proposition \ref{KernelpConv2} shows  that this vector tends to $a_2$ in $n$. Finally, we define the density $\check g^{(2)}_n$ as $\check g^{(2)}_n=g_n^{(1)}\frac{f_{\check a_2,n}}{g^{(1)}_{\check a_2,n}}$ which estimates $g^{(2)}$ through theorem \ref{KernelKRessultatPricipal}.\\
And so on, we will end up obtaining a sequence  $(\check a_1,\check a_2,...)$ of vectors in  $\R^d_*$ estimating the co-vectors of $f$ and a sequence of densities $(\check g^{(k)}_n)_k$ such that $\check g^{(k)}_n$ estimates $g^{(k)}$ through theorem \ref{KernelKRessultatPricipal}.
\vskip 3mm
\noindent {\bf3. Results}\label{OurResult}\\
\noindent {\bf3.1.}{\it  Convergence results}\label{cvResult}\\
\noindent {\bf3.1.1.}{\it  Hypotheses on $f$}\label{HypoF}

In this paragraph, we define the set of hypotheses on $f$ which can possibly be used in our work. Discussion on several of these hypotheses can be found in Appendix D.
In this section, to be more legible we replace $g$ with $g^{(k-1)}$. Let
$\Theta =\R^d_*$, $M(b,a,x)=\int ln(\frac{g(x)}{f(x)}\frac{f_b(\transp bx)}{g_b(\transp bx)})g(x)\frac{f_a(\transp ax)}{g_a(\transp ax)}dx-\ (\frac{g(x)}{f(x)}\frac{f_b(\transp bx)}{g_b(\transp bx)}-1)$,

$\Pn_n M(b,a)=\int M(b,a,x)d\Pn_n$, $\PP M(b,a)=\int M(b,a,x)f(x)dx,$\\
$\PP$ being the probability measure of $f$.
Similarly as in chapter $V$ of \cite{MR1652247}, we define :\\
$(H'1)$ : For all $\e>0$, there is $\eta>0$, such that for all $c\in\Theta$ verifying

$\ \ \ \ $ $\|c-a_k\|\geq \e$, we have $\PP M(c,a)<\PP M(a_k,a)-\eta$, with $a\in\Theta$.\\
$(H'2)$ : There exists a neighborhood of $a_k$, $V$, and a positive function $H$, such

$\ \ \ \ $that, for all $c\in V$ we have $|M(c,a_k,x)|\leq H(x)$ $(\PP -a.s.)$ with $\PP H<\infty$,\\
$(H'3)$ : There exists a neighborhood of $a_k$, $V$, such that for all $\e$, there exists a $\eta$  such that

$\ \ \ \ $ for all $c \in V$ and $a\in\Theta$, verifying $\|a-a_k\|\geq \e,$ we have $\PP M(c,a_k)<\PP M(c,a)-\eta$.\\
Putting $I_{a_k}=\frac{\dr^2}{\dr a^2}   K (g\frac{f_{a_k}}{g_{a_k}},f),$
and $x\to \rho(b,a,x)=ln(\frac{g(x)f_b(\transp bx)}{f(x)g_b(\transp bx)})\frac{g(x)f_a(\transp ax)}{g_a(\transp ax)}$, we now consider :\\
$(H'4)$ : There exists a neighborhood of $(a_k,a_k)$, $V_k'$, such that, for all  $(b,a)$ of $V_k'$, the gradient 

$\ \ \ \ $ $\nabla (\frac{g(x)f_a(\transp ax)}{g_a(\transp ax)})$ and the Hessian $\cH(\frac{g(x)f_a(\transp ax)}{g_a(\transp ax)})$ exist ($\lambda \_a.s.$), and the first order partial derivative

$\ \ \ \ $ $\frac{g(x)f_a(\transp ax)}{g_a(\transp ax)}$ and the first and second order derivative  of $(b,a)\mapsto \rho(b,a,x)$ are

$\ \ \ \ $ dominated  ($\lambda\_$a.s.) by integrable functions.\\
$(H'5)$ : The function $(b,a)\mapsto M(b,a,x)$ is $\cC^3$ in a neighborhood $V_k'$ of $(a_k,a_k)$ for all $x$ and all the  

$\ \ \ \ $partial derivatives of order 3 of $(b,a)\mapsto M(b,a,x)$ are dominated in $V_k'$ by a $\PP\_$integrable 

$\ \ \ \ $function $H(x)$.\\
$(H'6)$ : $\PP\|\frac{\dr}{\dr b}M(a_k,a_k)\|^2$ and $\PP\|\frac{\dr}{\dr a}M(a_k,a_k)\|^2$ are finite and the expressions

$\ \ \ \ $ $\PP\frac{\dr^2}{\dr b_i\dr b_j}M(a_k,a_k)$ and  $I_{a_k}$ exist and are invertible.\\
$(H'7)$ : There exists $k$ such that  ${\bf P} M(a_k,a_k)= 0$. \\
$(H'8)$ : $(Var_{{\bf P}}(M(a_k,a_k)))^{1/2}$ exists and is invertible.\\
$(H'0)$: $f$ and $g$ are assumed to be positive and bounded and such that $K(g,f)\geq\int|f(x)-g(x)|dx$.
\vskip 3mm
\noindent {\bf3.1.2.}{\it  Estimation of the first co-vector of $f$}\label{Estimofa1}

Let $\cR$ be the class of all positive functions $r$ defined on $\R$ and such that $g(x)r(\transp ax)$ is a density on $\R^d$ for all $a$ belonging to $\R^d_*$. The following proposition shows that there exists a vector $a$ such that $\frac{f_a}{g_a}$ minimizes $K(gr,f)$ in $r$:
\begin{proposition} \label{lemmeHuberModifprop}
There exists a vector $a$ belonging to $\R^d_*$ such that

$arg\min_{r\in\cR}K(gr,f)=\frac{f_a}{g_a}\text{ and }r(\transp ax)=\frac{f_a(\transp ax)}{g_a(\transp ax)}.$
\end{proposition}
Following \cite{MR2054155}, let us introduce the estimate of $K(g\frac{f_{a,n}}{g_{a}},f_{n})$, through 

$\check   K(g\frac{f_{a,n}}{g_{a}},f_{n})=\int M(a,a,x) d\Pn_n(x)$
\begin{proposition}\label{info}
Let $\check  a   :=  arg\inf_{a\in\R^d_*}\check   K(g\frac{f_{a,n}}{g_{a}},f_{n}).$\\
Then, $\check  a$ is a strongly convergent estimate of $a$, as defined in proposition \ref{lemmeHuberModifprop}.
\end{proposition}
\noindent Let us also introduce the following sequences  $(\check a_{k})_{k\geq 1}$ and $(\check g^{(k)}_{n})_{k\geq 1}$, for any given $n$ - see section 2.2.:\\
$\bullet$ $\check a_{k}$ is an estimate of $a_{k}$ as defined in proposition \ref{info}  with $\check g^{(k-1)}_{n}$ instead of $g$,\\
$\bullet$ $\check g^{(k)}_{n}$ is such that $\check g^{(0)}_{n}=g$, $\check g^{(k)}_{n}(x)=\check g^{(k-1)}_{n}(x)\frac{f_{\check a_k,n}(\transp {\check a_k}x)}{[\check g^{(k-1)}]_{\check a_k,n}(\transp {\check a_k}x)}$, i.e.
$\check g^{(k)}_{n}(x)=g(x)\Pi_{j=1}^k\frac{f_{\check a_j,n}(\transp {\check a_j}x)}{[\check g^{(j-1)}]_{\check a_j,n}(\transp {\check a_j}x)}$. We also note that $\check g^{(k)}_n$ is a density.
\vskip 3mm
\noindent {\bf3.1.3.}{\it  Convergence study at the  $k^{\text{th}}$ step of the algorithm:}

In this paragraph, we show that the sequence $(\check a_k)_n$ converges towards $a_k$ and that the sequence $(\check g^{(k)}_n)_n$ converges towards $g^{(k)}$.\\ 
Let $\check c_n(a)    =\ arg\sup_{c\in\Theta }\ \Pn_nM(c,a),$ with $a\in\Theta$,
and $\check \gamma_n   =\ arg\inf_{a\in\Theta }\ \sup_{c\in\Theta }\ \Pn_nM(c,a)$.
We state
\begin{proposition}\label{KernelpConv2}	Both $\sup_{a\in\Theta}\|\check c_n(a)-a_k\|$  and  $\check \gamma_n $ converge toward $a_k$ a.s.
\end{proposition}
\noindent Finally, the following theorem shows that $\check  g^{(k)}_n$ converges almost everywhere towards $g^{(k)}$:
\begin{theorem}\label{KernelKRessultatPricipal}It holds
$\check  g^{(k)}_n\to_n   g^{(k)}\ a.s.$
\end{theorem}
\noindent {\bf3.2.}{\it  Asymptotic inference at the  $k^{\text{th}}$ step of the algorithm}\label{AsymptResult}

The following theorem shows that $\check g^{(k)}_n$ converges towards $ g^{(k)}$ at the rate $O_{\PP}(m^{-\frac{1}{4+d}})$ in three differents cases, namely for any given $x$, with the  $L^1$ distance and with the Kullback-Leibler divergence:
\begin{theorem}\label{KernelSuperdiffconv}It holds 
$|\check g^{(k)}_n(x)-g^{(k)}(x)|=O_{\PP}(m^{-\frac{1}{4+d}}),$
$\int |\check g^{(k)}_n(x)-g^{(k)}(x)|dx=O_{\PP}(m^{-\frac{1}{4+d}})$ and
$|K(\check g^{(k)}_n,f)-K(g^{(k)},f)|=O_{\PP}(m^{-\frac{1}{4+d}}).$
\end{theorem}
Then, the following theorem shows that the laws of our estimators of $a_k$, namely $\check c_n(a_k)$ and $\check \gamma_n$, converge towards a linear combination of Gaussian variables.
\begin{theorem}\label{ KernelKloiestimateursMin}
It holds
$\sqrt n\cA.(\check c_n(a_k)-a_k)\cvL \cB.\cN_d(0,\PP\|\frac{\dr}{\dr b}M(a_k,a_k)\|^2)+\cC.\cN_d(0,\PP\|\frac{\dr}{\dr a}M(a_k,a_k)\|^2)$ and
$\sqrt n\cA.(\check \gamma_n-a_k)\cvL \cC.\cN_d(0,\PP\|\frac{\dr}{\dr b}M(a_k,a_k)\|^2)+\cC.\cN_d(0,\PP\|\frac{\dr}{\dr a}M(a_k,a_k)\|^2)$\\
where $\cA=\PP\frac{\dr^2}{\dr b\dr b}M(a_k,a_k)(\PP\frac{\dr^2}{\dr a\dr a}M(a_k,a_k)+\PP\frac{\dr^2}{\dr a\dr b}M(a_k,a_k))$, \\
$\cC=\PP\frac{\dr^2}{\dr b\dr b}M(a_k,a_k)$ and
$\cB=\PP\frac{\dr^2}{\dr b\dr b}M(a_k,a_k)+\PP\frac{\dr^2}{\dr a\dr a}M(a_k,a_k)+\PP\frac{\dr^2}{\dr a\dr b}M(a_k,a_k)$.
\end{theorem}
\noindent {\bf3.3.}{\it A stopping rule for the procedure}\label{StopAlgo}

In this paragraph, we show that $g_n^{(k)}$ converges towards $f$ in $k$ and $n$. Then, we provide a stopping rule for this identification procedure.
\vskip 3mm
\noindent {\bf3.3.1.}{\it Estimation of $f$}\label{EstOfF}

Through remark \ref{criteria-H} and as explained in section 14 of \cite{MR790553}, the following lemma shows that $K(g^{(k-1)}_n\frac{f_{a_k,n}}{g^{(k-1)}_{a_k,n}},f_{a_k,n})$ converges almost everywhere towards zero as $k$ goes to infinity and thereafter as $n$ goes to infinity : 
\begin{lemma}\label{FromSection14Huber}
We have $\lim_n\lim_kK(\check g^{(k)}_n\frac{f_{{a_k},n}}{[\check g^{(k)}]_{{a_k},n}},f_n)= 0$ a.s.
\end{lemma}
\noindent Consequently, the following proposition provides us with an estimate of $f$:
\begin{theorem}\label{limnk}
We have $\lim_n\lim_k\check g^{(k)}_n=f$ a.s.
\end{theorem}
\noindent {\bf3.3.2.}{\it Testing of the criteria}\label{Test}

In this paragraph, through a test of the criteria,  namely $a\mapsto K(\check g^{(k)}_n\frac{f_{a,n}}{[\check g^{(k)}]_{a,n}},f_n)$, we build a stopping rule for this identification procedure. First, the next theorem enables us to derive the law of the criteria:
\begin{theorem} \label{KernelLOIDUCRITERE}
For a fixed $k$, we have 

$\sqrt n(Var_{\PP}(M(\check c_n(\check \gamma_n),\check \gamma_n)))^{-1/2}(\Pn_nM(\check c_n(\check \gamma_n),\check \gamma_n)-\Pn_nM(a_k,a_k)) \cvL \cN(0,I)$,\\
as $n$ goes to infinity, where $k$ represents the $k^{th}$ step of the algorithm and $I$ is the identity matrix in $\R^d$.
\end{theorem}
Note that $k$ is fixed in theorem \ref{KernelLOIDUCRITERE} since $\check \gamma_n   =\ arg\inf_{a\in\Theta }\ \sup_{c\in\Theta }\ \Pn_nM(c,a)$ where $M$ is a known function of $k$, see section 3.1.1. Thus, in the case where $K(g^{(k-1)}\frac{f_{a_k}}{g^{(k-1)}_{a_k}},f)= 0$, we obtain
\begin{corollary} \label{KernelLOIDUCRITERE2} 
We have 
$\sqrt n(Var_{\PP}(M(\check c_n(\check \gamma_n),\check \gamma_n)))^{-1/2}(\Pn_nM(\check c_n(\check \gamma_n),\check \gamma_n)) \cvL \cN(0,I)$.
\end{corollary}
Hence, we propose the test of the null hypothesis 

$\text{$(H_0)$ : $K(g^{(k-1)}\frac{f_{a_k}}{g^{(k-1)}_{a_k}},f)= 0$ versus $(H_1)$ : $K(g^{(k-1)}\frac{f_{a_k}}{g^{(k-1)}_{a_k}},f)\not= 0$.}$\\
Based on this result, we stop the algorithm, then, defining $a_k$ as the last vector generated, we derive from corollary \ref{KernelLOIDUCRITERE2} a $\alpha$-level confidence ellipsoid around $a_k$, namely 

$\cE_k=\{b\in\R^d;\ \sqrt n(Var_{\PP}(M(b,b)))^{-1/2}\Pn_nM(b,b)\leq q_{\alpha}^{\cN(0,1)} \}$\\
where $q_{\alpha}^{\cN(0,1)}$ is the quantile of a $\alpha$-level reduced centered normal distribution and where $\Pn_n$ is the empirical measure araising from a realization of the sequences $(X_1,\ldots,X_n)$ and $(Y_1,\ldots,Y_n)$.\\ The following corollary thus provides us with a confidence region for the above test:
\begin{corollary}\label{KernelLOIDUCRITERE2coro}
$\cE_k$ is a confidence region for the test of the null hypothesis $(H_0)$ versus  $(H_1)$.
\end{corollary}
\noindent {\bf4. Comparison of all the optimisation methods}

In this section, we study Huber's algorithm in a similar manner to sections 2 and 3. We will then be able to compare our methodologies.\\
Until now, the choice has only been to use the class of Gaussian distributions. Here and similarly to section 2.1, we extend this choice to the class of elliptical distributions.
Moreover, using the subsample $X_1$, $X_2$,..., $X_n$, see Appendix B, and using  the procedure of section 2.2. with $K(g_a,f_a)$, see section 4.2, instead of $K(g\frac{g_a}{f_a},f)$, proposition \ref{QuotientDonneLoi}, lemma \ref{toattaint} and remark \ref{criteria-H} enable us to perform the Huber's algorithm :\\
$\bullet$ we define $\hat a_1$ and the density $\hat g^{(1)}_n$ such that $\hat a_1=arg\max_{a\in\R^d_*}K(g_{a},f_{a,n})$ and $\hat g^{(1)}_n=g\frac{f_{\hat a_1,n}}{g_{\hat a_1}}$,\\
$\bullet$ we define $\hat a_2$ and the density $\hat g^{(2)}_n$ such that $\hat a_2=arg\max_{a\in\R^d_*}K(\hat g^{(1)}_{a,n},f_{a,n})$ and $\hat g^{(2)}_n=\hat g^{(1)}_n\frac{f_{\hat a_2,n}}{\hat g^{(1)}_{\hat a_2,n}}$,\\
and so on, we obtain a sequence  $(\hat a_1,\hat a_2,...)$ of vectors in  $\R^d_*$  and a sequence of densities $\hat g^{(k)}_n$.
\vskip 3mm
\noindent {\bf4.1.}{\it  Hypotheses on $f$}\label{HypoFHuber}

In this paragraph, we define the set of hypotheses on $f$ which can be of use in our present work.
First, we denote $g$ in lieu of $g^{(k-1)}$. Let
$\Theta_a^1 =\{b\in\Theta\ |\ \ \int\ (\frac{g_b(\transp bx)}{f_b(\transp bx)}-1)f_a(\transp ax)\ dx<\infty\},$

$m(b,a,x)=\int\ ln(\frac{g_b(\transp bx)}{f_b(\transp bx)})g_a(\transp ax)\ dx\ -\ (\frac{g_b(\transp bx)}{f_b(\transp bx)}-1)$,

$\PP^a m(b,a)=\int\ m(b,a,x)f_a(\transp ax)\ dx$ and $\Pn_n m(b,a)=\int\ m(b,a,x) \frac{f_a(\transp ax)}{f(x)}d\Pn_n$,\\
$\PP^a $ being the probability measure of $f_a$.\\
Similarly as in chapter $V$ of  \cite{MR1652247}, we define :\\
$(H1)$ : For all $\e>0$, there is $\eta>0$ such that, for all $b\in\Theta_a^1$ verifying

$\ \ \ \ \ $ $\|b-a_k\|\geq \e$ for all $a\in\Theta$, we have $\PP^a m(b,a)<\PP^a m(a_k,a)-\eta$,\\
$(H2)$ : There exists a neighborhood of $a_k$, $V$, and a positive function $H$, such

$\ \ \ $ that, for all $b\in V$, we have $|m(b,a_k,x)|\leq H(x)$ $(\PP^a -a.s.)$ with $\PP^a H<\infty$,\\
$(H3)$ : There exists a neighborhood $V$ of $a_k$, such that for all $\e$, there exists a $\eta$  such

$\ \ \ \ $ that for all $b \in V$  and $a\in\Theta$, verifying $\|a-a_k\|\geq \e$, we have $\PP^{a_k}m(b,a_k)-\eta>\PP^am(b,a).$\\
Moreover, defining  $x\to \upsilon(b,a,x)=ln(\frac{g_b(\transp bx)}{f_b(\transp bx)})g_a(\transp ax)$, putting:\\
$(H4)$ : There exists a neighborhood  of $(a_k,a_k)$, $V_k$, such that, for all $(b,a)$ of $V_k$,

$\ \ \ \ $ the gradient $\nabla (\frac{g_a(\transp ax)}{f_a(\transp ax)})$ and the Hessian $\cH(\frac{g_a(\transp ax)}{f_a(\transp ax)})$ exist ($\lambda-a.s.$) and the first order  partial

$\ \ \ \ $  derivative $\frac{g_a(\transp ax)}{f_a(\transp ax)}$ and the first and second order derivative  of order 3 of $(b,a)\mapsto \upsilon(b,a,x)$ 

$\ \ \ $ are dominated ($\lambda\_$a.s.) by integrable functions.\\
$(H5)$ : The function $(b,a)\mapsto m(b,a)$ is $\cC^3$ in a neighborhood $V_k$ of $(a_k,a_k)$ for all $x$ and all the 

$\ \ \ \ $ partial derivatives of $(b,a)\mapsto m(b,a)$ are dominated in $V_k$ by a $\PP\_$integrable function $H(x)$.\\
$(H6)$ : $\PP\|\frac{\dr}{\dr b}m(a_k,a_k)\|^2$ and $\PP\|\frac{\dr}{\dr a}m(a_k,a_k)\|^2$ are finite and the quantities 

$\ \ \ $ $\PP\frac{\dr^2}{\dr b_i\dr b_j}m(a_k,a_k)$ and $\PP\frac{\dr^2}{\dr a_i\dr a_j}m(a_k,a_k)$ are invertible.\\
$(H7)$ : There exists $k$ such that ${\bf P} m(a_k,a_k)= 0$.\\
$(H8)$ : $(Var_{\PP}(m(a_k,a_k)))^{1/2}$ exists and is invertible.
\vskip 3mm
\noindent {\bf4.2.}{\it  The first co-vector of $f$ simultaneously optimizes four problems}\label{a1solve4}

We first study Huber's analytic approach. Let $\cR'$ be the class of all positive functions $r$ defined on $\R$ and such that $f(x)r^{-1}(\transp ax)$ is a density on $\R^d$  for all $a$ belonging to $\R^d_*$. The following proposition shows that there exists a vector $a$ such that  $\frac{f_a}{g_a}$ minimizes $K(fr^{-1},g)$ in $r$:
\begin{proposition}[Analytic Approach]\label{lemmeHuber0prop}
There exists a vector $a$ belonging to $\R^d_*$ such that 

$arg\min_{r\in\cR'}K(fr^{-1},g)=\frac{f_a}{g_a},\ r(\transp ax)=\frac{f_a(\transp ax)}{g_a(\transp ax)}$ as well as $  K(f,g)=  K(f_a,g_a)+  K(f\frac{g_a}{f_a},g).$
\end{proposition}
\noindent We also study Huber's synthetic approach. Let $\cR$ be the class of all positive functions $r$ defined on $\R$ and such that $g(x)r(\transp ax)$ is a density on $\R^d$ for all $a$ belonging to $\R^d_*$. The following proposition shows that there exists a vector $a$ such that  $\frac{f_a}{g_a}$ minimizes $K(gr,f)$ in $r$:
\begin{proposition}[Synthetic Approach]\label{lemmeHuberprop}
There exists a vector $a$ belonging to $\R^d_*$ such that

$arg\min_{r\in\cR}K(f,gr)=\frac{f_a}{g_a},\  r(\transp ax)=\frac{f_a(\transp ax)}{g_a(\transp ax)}$
 as well as  $  K(f,g)=  K(f_a,g_a)+  K(f,g\frac{f_a}{g_a}).$
\end{proposition}
\noindent In the meanwhile, the following proposition shows  that there exists a vector $a$ such that   $\frac{f_a}{g_a}$ minimizes  $K(g,fr^{-1})$ in $r$.
\begin{proposition} \label{liendeminimaxi}
There exists a vector $a$ belonging to $\R^d_*$ such that
$arg\min_{r\in\cR'}K(g,fr^{-1})=\frac{f_a}{g_a},$  and $r(\transp ax)=\frac{f_a(\transp ax)}{g_a(\transp ax)}$.
Moreover, we  have $K(g,f)=  K(g_a,f_a)+K(g,f\frac{g_a}{f_a})$.
\end{proposition}
\begin{remark}\label{criteria-H}First, through property \ref{ExitenceDeLEntropieDesProj}, we get $K(f,g\frac{f_a}{g_a})= K(g,f\frac{g_a}{f_a})= K(f\frac{g_a}{f_a},g)$ and $K(f_a,g_a)=K(g_a,f_a)$. Thus, proposition \ref{liendeminimaxi} implies that finding the argument of the maximum of $K(g_a,f_a)$ amounts to finding the argument of the maximum of $K(f_a,g_a)$. Consequently, the criteria of Huber's methodologies is $a\mapsto K(g_a,f_a)$.
Second, our criteria is $a\mapsto K(g\frac{g_a}{f_a},f)$ and property \ref{ExitenceDeLEntropieDesProj} implies $K(g,f\frac{g_a}{f_a})= K(g\frac{f_a}{g_a},f)$.
Consequently, since \cite{MR2054155} takes into account the very form of the criteria, we are then in a position to compare Huber's methodologies with ours.
\end{remark}
To recapitulate, the choice of $r=\frac{f_a}{g_a}$ enables us to simultaneously solve the following four optimisation problems, for $a\in\R^d_*$:

$\text{First, find $a$ such that }a\ =\ {\text{\itshape arginf}}_{a\in\R^d_*}\ K(f\frac{g_a}{f_a},g)$ 
- analytic approach -

$\text{Second, find $a$ such that }a\ =\ {\text{\itshape arginf}}_{a\in\R^d_*}\ K(f,g\frac{f_a}{g_a})$ - synthetic approach -

$\text{Third, find $a$ such that }a\ =\ {\text{\itshape argsup}}_{a\in\R^d_*}\ K(g_a,f_a)$ 
- to compare Huber's methods with ours -

$\text{Fourth, find $a$ such that }a\ =\ {\text{\itshape arginf}}_{a\in\R^d_*}\ K(g\frac{f_a}{g_a},f)$ - our method.
\vskip 3mm
\noindent {\bf4.2.}{\it  On the sequence of the transformed densities $(g^{(j)})$}

As already explained in the introduction section, the Mu Zhu article leads us to only consider Huber's synthetic approach.

\noindent {\bf4.2.1.}{\it Estimation of the first co-vector of $f$}

Using the subsample $X_1$, $X_2$,..,$X_n$, see Appendix B, and following \cite{MR2054155}, let us introduce the estimate of $K(g_{a},f_{a,n})$, through $\hat   K(g_{a},f_{a,n})=\int m(a,a,x)(\frac{f_{a,n}(\transp ax)} {f_n(x)}) d\Pn_n$
\begin{proposition}\label{intox}
Let $\hat a := arg\sup_{a\in\R^d_*}\hat   K(g_{a},f_{a,n}).$\\ Then, $\hat a$ is a strongly convergent estimate of $a$, as defined in proposition \ref{liendeminimaxi}.
\end{proposition}
\noindent Finally, we define the following sequences  $(\hat a_{k})_{k\geq 1}$ and $(\hat g^{(k)}_n)_{k\geq 1}$ - for any given $n$ :\\
$\bullet$ $\hat a_{k}$ is an estimate of $a_{k}$ as defined in proposition $\ref{intox}$ with $\hat g^{(k-1)}_n$ instead of $g$,\\
$\bullet$ $\hat g^{(k)}_n$ is such that $\hat g^{(0)}_n=g$ and $\hat g^{(k)}_n(x)=\hat g^{(k-1)}_n(x)\frac{f_{\hat a_k,n}(\transp {\hat a_k}x)}{[\hat g^{(k-1)}]_{\hat a_k,n}(\transp {\hat a_k}x)}$, i.e.
$\hat g^{(k)}_n(x)=g(x)\Pi_{j=1}^k\frac{f_{\hat a_j,n}(\transp {\hat a_j}x)}{[\hat g^{(j-1)}]_{\hat a_j,n}(\transp {\hat a_j}x)}$.
\vskip 3mm
\noindent {\bf4.2.2.}{\it  Convergence study at the  $k^{\text{th}}$ step of the algorithm}

Let $\hat b_n(a)=\ arg\sup_{b\in\Theta }\ \Pn^a_nm(b,a)$, with $a\in\Theta$, and $\hat \beta_n=\ arg\sup_{a\in\Theta }\ \sup_{b\in\Theta }\ \Pn^a_nm(b,a)$, then
\begin{proposition}\label{pConv1}
Both $\sup_{a\in\Theta}\|\hat b_n(a)-a_k\| $  and  $\hat \beta_n $ converge toward $a_k$ a.s.
\end{proposition}
\noindent Finally, the following theorem shows that  $\hat  g^{(k)}_n$ converges almost everywhere towards $g^{(k)}$ :
\begin{theorem}\label{KernelKRessultatPricipal-H}
For any given $k$, it holds $\hat  g^{(k)}_n\to_n   g^{(k)}$ a.s.
\end{theorem}
\noindent {\bf4.2.3.}{\it Asymptotic inference at the  $k^{\text{th}}$ step of the algorithm}\label{AsymptResultHuber}

The following theorem shows that $\hat g^{(k)}_n$ converges towards $g^{(k)}$ at the rate $O_{\PP}(m^{-\frac{1}{4+d}})$ in three differents cases, namely for any given $x$, with the  $L^1$ distance and with the Kullback-Leibler divergence:
\begin{theorem}\label{KernelSuperdiffconv-H}It holds $|\hat g^{(k)}_n(x)-g^{(k)}(x)|=O_{\PP}(m^{-\frac{1}{4+d}})$, 
$\int |\hat g^{(k)}_n(x)-g^{(k)}(x)|dx=O_{\PP}(m^{-\frac{1}{4+d}})$ and 
$|K(f,\hat g^{(k)}_n)-K(f,g^{(k)})|=O_{\PP}(m^{-\frac{1}{4+d}})$.
\end{theorem}
The following theorem shows that the laws of Huber's estimators of $a_k$, namely $\hat b_n(a_k)$ and $\hat \beta_n$, converge towards a linear combination of Gaussian variables.
\begin{theorem} \label{ KloiestimateursMax}
It holds
$\sqrt n\cD.(\hat b_n(a_k)-a_k) \cvL \cE.\cN_d(0,\PP\|\frac{\dr}{\dr b}m(a_k,a_k)\|^2)+\cF.\cN_d(0,\PP\|\frac{\dr}{\dr a}m(a_k,a_k)\|^2)$ and
$\sqrt n\cD.(\hat \beta_n-a_k)  \cvL \cG.\cN_d(0,\PP\|\frac{\dr}{\dr a}m(a_k,a_k)\|^2)+\cF.\cN_d(0,\PP\|\frac{\dr}{\dr b}m(a_k,a_k)\|^2)$\\
where $\cE=\PP\frac{\dr^2}{\dr a^2}m(a_k,a_k)$, $\cF=\PP\frac{\dr^2}{\dr a\dr b}m(a_k,a_k),$ $\cG=\PP\frac{\dr^2}{\dr b^2}m(a_k,a_k)$ and\\
$\cD=(\PP\frac{\dr^2}{\dr b^2}m(a_k,a_k)\PP\frac{\dr^2}{\dr a^2}m(a_k,a_k)-\PP\frac{\dr^2}{\dr a\dr b}m(a_k,a_k)\PP\frac{\dr^2}{\dr b\dr a}m(a_k,a_k))>0$.
\end{theorem}
\noindent {\bf4.3.}{\it A stopping rule for the procedure}

We first give an estimate of $f$. Then, we provide a stopping rule for this identification procedure.
\begin{remark}
In the case where $f$ is known, as explained in section 14 of \cite{MR790553}, the sequence $(K(g^{(k-1)}_{a_k},f_{a_k}))_{k\geq 1}$ converges towards zero. Many authors have studied this hypothesis and its consequences. For example, Huber deducts that, if $f$ can be deconvoluted with a Gaussian component,\\ $(K(g^{(k-1)}_{a_k},f_{a_k}))_{k\geq 1}$ converges toward $0$. He then shows that $g^{(i)}$ uniformly converges in $L^1$ towards $f$ - see propositions 14.2 and 14.3 page 461 of his article.
\end{remark}
\noindent {\bf4.3.1.}{\it  Estimation of $f$}

The following lemma shows that $\lim_k K(\hat g^{(k)}_{{a_k},n},f_{{a_k},n})$ converges towards zero as $k$ goes to infinity and thereafter as $n$ goes to infinity :
\begin{lemma}\label{KernelHuber0}
We have
$\lim_n\lim_k K(\hat g^{(k)}_{{a_k},n},f_{{a_k},n})= 0,$ a.s.
\end{lemma}
\noindent Then, the following theorem enables us to provide simulations through an estimation of $f$ 
\begin{theorem}\label{limnk-H}
We have $\lim_n\lim_k\hat g^{(k)}_n=f,$  a.s.
\end{theorem}
\noindent {\bf4.3.2.}{\it  Testing of the criteria}

In this paragraph, through a test of Huber's criteria, namely $a\mapsto K(\hat g^{(k)}_{a,n},f_{a,n})$, we will build a stopping rule for the procedure. First, the next theorem gives us the law of Huber's criteria.
\begin{theorem} \label{KernelLOIDUCRITERE-H} For a fixed $k$, we have

$\sqrt n(Var_{\PP}(m(\hat b_n(\hat \beta_n),\hat \beta_n)))^{-1/2}(\Pn_nm(\hat b_n(\hat \beta_n),\hat \beta_n)-\Pn_nm(a_k,a_k)) \cvL \cN(0,I)$, \\
as $n$ goes to infinity, where $k$ represents the $k^{th}$ step of the algorithm and $I$ is the identity matrix in $\R^d$.
\end{theorem}
\noindent Note that $k$ is fixed in theorem \ref{KernelLOIDUCRITERE-H} since $\hat \beta_n=\ arg\sup_{a\in\Theta }\ \sup_{b\in\Theta }\ \Pn^a_nm(b,a)$ where $m$ is a known function of $k$ - see section 4.1. Thus, in the case where $K(g^{(k)}_{a},f_{a})=0$, we obtain 
\begin{corollary} \label{LOIDUCRITERE2-H} 
$\\$We have  $\sqrt n(Var_{\PP}(m(\hat b_n(\hat \beta_n),\hat \beta_n)))^{-1/2}(\Pn_nm(\hat b_n(\hat \beta_n),\hat \beta_n)) \cvL \cN(0,I).$
\end{corollary}
\noindent Hence, we propose the test of the null hypothesis $(H_0)$ : $K(g^{(k-1)}_{a_k},f_{a_k})= 0$ versus the alternative $(H_1)$ : $K(g^{(k-1)}_{a_k},f_{a_k})\not= 0$. Based on this result, we stop the algorithm, then, defining $a_k$ as the last vector generated from the Huber's algorithm, we derive from corollary \ref{LOIDUCRITERE2-H}, a $\alpha$-level confidence ellipsoid around $a_k$, namely $\cE_k'=\{b\in\R^d;\ \sqrt n(Var_{\PP}(m(b,b)))^{-1/2}\Pn_nm(b,b)\leq q_{\alpha}^{\cN(0,1)} \}$
where $q_{\alpha}^{\cN(0,1)}$ is the quantile of a $\alpha$-level reduced centered normal distribution and where $\Pn_n$ is the empirical measure araising from a realization of the sequences $(X_1,\ldots,X_n)$ and $(Y_1,\ldots,Y_n)$.\\ Consequently, the following corollary provides us with a confidence region for the above test:
\begin{corollary}\label{LOIDUCRITERE2coro-H}
$\cE_k'$ is a confidence region for the test of the null hypothesis $(H_0)$ versus $(H_1)$.
\end{corollary}
\noindent {\bf5. Simulations}\label{Simul400}

We  illustrate this section by detailing three simulations.\\
In each simulation, the program  follows our algorithm and aims at creating a sequence of densities $(g^{(j)})$, $j=1,..,k$, $k<d$, such that
$g{(0)}=g,$ $g^{(j)}=g^{(j-1)}f_{a_j}/[g^{(j-1)}]_{a_j}$ and $K(g^{(k)},f)=0,$
where $K$ is the Kullback-Leibler divergence and $a_j=arg\inf_b K(g^{(j-1)}f_b/[g^{(j-1)}]_b,f)$, for all $j=1,...,k$.\\
Then, in the first two simulations, the  program follows Huber's method and generates a sequence of densities $(g^{(j)})$, $j=1,..,k$, $k<d$, such that $g{(0)}=g,\ g^{(j)}=g^{(j-1)}f_{a_j}/[g^{(j-1)}]_{a_j}$ and $K(f,g^{(k)})=0,$
where $K$ is the Kullback-Leibler divergence and  $a_j=argsup_b K([g^{(j-1)}]_b,f_b)$, for all $j=1,...,k$.\\ Finally, in the third example, we study the robustness of our method with four outliers.
\begin{Sim} \label{Sim1}
\em$\\$ We are in dimension $3$(=d). We consider a sample of $50$(=n) values of a random variable $X$ with density $f$ defined by, 

$f(x)=Normal(x_1+x_2).Gumbel(x_0+x_2).Gumbel(x_0+x_1)$,\\
where the Gumbel law parameters are $(-3, 4)$ and $(1,1)$ and where the normal distribution parameters are $(-5,2)$.We generate a Gaussian random variable $Y$ with a density - that we will name $g$ - which has the same mean and variance as $f$.\\
In the first part of the  program, we theoretically obtain $k=2$, $a_1=(1,0,1)$ and $a_2=(1,1,0)$ (or $a_2=(1,0,1)$ and $a_1=(1,1,0)$ which leads us to the same conclusion). To get  this result, we perform the following test 

$(H_0):\ (a1,a2)=((1,0,1),(1,1,0))\ versus\ (H_1):(a1,a2)\not=((1,0,1),(1,1,0)).$\\
Moreover, if $i$ represents the last iteration of the algorithm, then 

$\sqrt n (Var_\PP(M(c_n(\gamma_n),\gamma_n)))^{(-1/2)}\Pn_nM(c_n(\gamma_n),\gamma_n)\cvL \cN(0,1),$\\
and then we estimate $(a_1,a_2)$ with the following 0.9(=$\alpha$) level confidence ellipsoid

$\cE_i=\{b \in \R^3;\  (Var_\PP(M(b,b)))^{-1/2}\Pn_nM(b,b)\leq q^{\cN(0,1)}_{\alpha}/\sqrt n\simeq\frac{0,2533}{7.0710678}=0.03582203\}.$\\
Indeed, if $i=1$ represents the last iteration of the algorithm, then $a_1\in \cE_{1}$, and if $i=2$ represents the last iteration of the algorithm, then $a_2\in \cE_{2}$, and so on, if $i$ represents the last iteration of the algorithm, then $a_i\in \cE_{i}$.

Now, if we follow Huber's method, we also theoretically obtain $k=2$, $a_1=(1,0,1)$ and $a_2=(1,1,0)$ (or $a_2=(1,0,1)$ and $a_1=(1,1,0)$ which leads us to the same conclusion).
To get this result, we perform the following test:

$(H_0):\ (a_1,a_2)=((1,0,1),(1,1,0))\ versus\ (H_1):\ (a_1,a_2)\not=((1,0,1),(1,1,0)).$\\
Similarly as above, the fact that, if $i$ represents the last iteration of the algorithm, then

$\sqrt n(Var_\PP(m(b_n(\beta_n),\beta_n)))^{(-1/2)}\Pn_nm(b_n(\beta_n),\beta_n)\cvL \cN(0,1),$
enables us to estimate our sequence of $(a_i)$, reduced to $(a_1,a_2)$, through the following 0.9(=$\alpha$) level confidence ellipsoid

$\cE_i'=\{b \in \R^3;\  (Var_\PP(m(b,b)))^{-1/2}\Pn_nm(b,b)\leq q^{\cN(0,1)}_{\alpha}/\sqrt n\simeq 0.03582203\}.$\\
Finally, we obtain
\begin{table}[h!]
\caption{Simulation 1 : Numerical results  of the optimisation.}
\label{Sim1-tab}
\begin{tabular}{cll}
\hline\noalign{\smallskip}
& Our Algorithm & Huber's Algorithm\\
\noalign{\smallskip}\hline\noalign{\smallskip}
\multirow{3}{8.5cm}{Projection Study 0 :}
 & minimum : 0.317505     & maximum : 0.715135 \\  \cline{2-3}
 & at point : (1.0,1.0,0)& at point : (1.0,1.0,0) \\  \cline{2-3}
 & P-Value : 0.99851       & P-Value : 0.999839 \\
\noalign{\smallskip}\hline
\multirow{1}{8.5cm}{Test :}
 & $H_0$ : $a_1\in \cE_{1}$ : False   & $H_0$ : $a_1\in \cE_{1}'$ : False \\
\noalign{\smallskip}\hline
\multirow{3}{8.5cm}{Projection Study 1 :}
 & minimum : 0.0266514     & maximum : 0.007277 \\  \cline{2-3}
 & at point : (1.0,0,1.0)       & at point : (1,0.0,1.0) \\  \cline{2-3}
 & P-Value : 0.998852      & P-Value : 0.999835 \\
\noalign{\smallskip}\hline
\multirow{1}{8.5cm}{Test :}
 & $H_0$ : $a_2\in \cE_{2}$ : True   & $H_0$ : $a_2\in \cE_{2}'$ : True \\
\noalign{\smallskip}\hline
\multirow{1}{8.5cm}{K(Estimate $g_m^{(2)}$, $g^{(2)}$)}
&0.444388&0.794124\\
\noalign{\smallskip}\hline
\end{tabular}
$\\$Therefore, we conclude that $f=g^{(2)}.$
\end{table}
\end{Sim}
\begin{Sim}\label{Sim2}
\em$\\$We are in dimension $10$(=d). We consider a sample of $50$(=n) values of a random variable $X$ with density $f$ defined by, 

$f(x)=Gumbel(x_0).Normal(x_1,\ldots,x_9)$,\\
where the Gumbel law parameters are -5 and 1 and where the  normal distribution is reduced and centered.\\
Our reasoning is the same as in Example \ref{Sim1}. In the first part of the  program, we theoretically obtain $k=1$ and $a_1=(1,0,\ldots,0)$. To get  this result, we perform the following test 

$(H_0):\ a_1=(1,0,\ldots,0)\ versus\ (H_1):\ a_1\not=(1,0,\ldots,0)$.\\
We estimate $a_1$ by the following 0.9(=$\alpha$) level confidence ellipsoid

$\cE_{i}=\{b \in \R^2;\  (Var_\PP(M(b,b)))^{-1/2}\Pn_nM(b,b)\leq q^{\cN(0,1)}_{\alpha}/\sqrt n\simeq 0.03582203\}.$\\
Now, if we follow Huber's method, we also theoretically obtain $k=1$ and $a_1=(1,0,\ldots,0)$.
To get this result, we perform the following test 

$(H_0):\ a_1=(1,0,\ldots,0)\ versus\ (H_1):\ a_1\not=(1,0,\ldots,0).$\\
Hence, using the same reasoning as in Example \ref{Sim1}, we  estimate $a_1$ through the following 0.9 (=$\alpha$) level confidence ellipsoid

$\cE_{i}'=\{b \in \R^2;\  (Var_\PP(m(b,b)))^{-1/2}\Pn_nm(b,b)\leq q^{\cN(0,1)}_{\alpha}/\sqrt n\simeq 0.03582203\}$.\\
And, we obtain
\begin{table}[h!]
\caption{Simulation 2 : Numerical results  of the optimisation.}
\label{Sim2-tab}
\begin{tabular}{cll}
\hline\noalign{\smallskip}
& Our Algorithm & Huber's Algorithm\\
\noalign{\smallskip}\hline\noalign{\smallskip}
\multirow{8}{5cm}{Projection Study  0:}
 & minimum : 0.00263554     & maximum : 0.00376235 \\  \cline{2-3}
 & at point : (1.0001, & at point : (0.9902, \\ 
 &         0.0040338, 0.098606, 0.115214, & 0.0946806, 0.161447, 0.0090245,\\
 &          0.067628, 0.16229, 0.00549203, &  0.147804, 0.180259, 0.0975065,\\
 &         0.014319, 0.149339, 0.0578906)  & 0.101044, 0.190976, 0.155706)\\  \cline{2-3}
 & P-Value : 0.828683       & P-Value : 0.807121 \\
\noalign{\smallskip}\hline
\multirow{1}{5cm}{Test :}
 & $H_0$ : $a_1\in \cE_{1}$ : True   & $H_0$ : $a_1\in \cE_{1}'$ : True \\
\noalign{\smallskip}\hline
\multirow{1}{5cm}{K(Estimate $g_m^{(1)}$, $g^{(1)}$)}
& 2.44546 & 2.32331\\
\noalign{\smallskip}\hline
\end{tabular}
$\\$Therefore, we conclude that $f=g^{(1)}.$
\end{table}
\end{Sim}
\begin{Sim}\label{Sim3}
\em$\\$We are in dimension $20$(=d). We first generate a sample with
$100$(=n) observations, namely four outliers $x=(2,0,\ldots,0)$ and
$96$ values of a random variable $X$ with a density $f$ defined by

$f(x)=Gumbel(x_0).Normal(x_1,\ldots,x_{19})$\\
where the Gumbel law parameters are -5 and 1 and where the normal
distribution is reduced and centered. Our reasoning is the same as
in Simulation \ref{Sim1}. \\We
theoretically obtain $k=1$ and $a_1=(1,0,\ldots,0)$. To get  this
result, we perform the following test

$(H_0):\ a_1=(1,0,\ldots,0)\ versus\ (H_1):\ a_1\not=(1,0,\ldots,0)$ \\
We estimate $a_1$ by the following 0.9(=$\alpha$) level confidence ellipsoid

$\cE_{i}=\{b \in \R^2;\  (Var_\PP(M(b,b)))^{-1/2}\Pn_nM(b,b)\leq q^{\cN(0,1)}_{\alpha}/\sqrt n\simeq 0.02533\}$\\ And, we obtain
\begin{table}[h]
\caption{Simulation 3: Numerical results of the optimisation.}
\label{Sim3-opti}
\begin{tabular}{cll}
\hline\noalign{\smallskip}
Our Algorithm &&\\
\noalign{\smallskip}\hline\noalign{\smallskip}
\multirow{5}{6cm}{Projection Study 0}
 & minimum  : 0.024110                                               \\  \cline{2-3}
 & at point : (0.8221, 0.0901, 0.0892, -0.2020, 0.0039, 0.1001,      \\
 & 0.0391, 0.08001, 0.07633, -0.0437, 0.12093, 0.09834, 0.1045,      \\
 & 0.0874, -0.02349, 0.03001,  0.12543, 0.09435, 0.0587, -0.0055)    \\ \cline{2-3}
 & P-Value  : 0.77004                                                \\
\noalign{\smallskip}\hline
\multirow{1}{6cm}{Test :}
 & $H_0$    : $a_1\in \cE_{1}$ : True                                \\
\noalign{\smallskip}\hline
\multirow{1}{6cm}{K(Estimate $g_m^{(1)}$, $g^{(1)}$)}
& 2.677015                                                           \\
\noalign{\smallskip}\hline
\end{tabular}
$\\$Therefore, we conclude that $f=g^{(1)}.$
\end{table}
\end{Sim}
\noindent{\bf Critics of the simulations}\\
As customary in simulation studies, as approximations accumulate, results depend on the power of the calculators used as well as on the available memory.
Moreover, in order to implement our optimisation in $\R^d$ of the relative entropy, we choose to apply the simulated annealing method.\\
Thus, in the case where $f$ is unknown, we will never have the certainty to have reached the desired minimum or maximum of the Kullback-Leibler divergence.
Indeed, this probabilistic metaheuristic only converges, and the probability to reach the minimum or the maximum only tends towards 1, when the number of random jumps tends in theory towards infinity.\\
We also note that no theory on the optimal number of jumps to implement does exist, as this number depends on the specificities of each particular problem.\\
Finally, we choose the $50^{-\frac{4}{4+d}}$ (resp. $100^{-\frac{4}{4+d}}$) for the AMISE of the simulations 1 and 2 (resp. 3). This choice leads us to simulate 50 (resp.100) random variables, see \cite{MR1191168} page 151, none of which have been discarded to obtain the truncated sample.

\noindent {\bf Conclusion}\\
Characteristic structures as well as one-dimensional projections and their associated distributions in multivariate datasets can be evidenced through Projection Pursuit.\\
The present article demonstrates that our Kullback-Leibler divergence minimisation method constitutes a good alternative to Huber's relative entropy maximization approach, see \cite{MR790553}. Indeed, the convergence results as well as the simulations we carried out clearly evidences the robustness of our methodology.
\vskip 3mm
\noindent {\bf A. Reminders}\\
\noindent {\bf A.1.}{\it The relative entropy (or Kullback-Leibler divergence)}

We call $h_a$ the density of  $\transp a Z$ if  $h$ is the density of $Z$, and $K$ the relative entropy or Kullback-Leibler divergence. The function $K$ is defined by - considering $P$ and $Q$, two probabilities:

$K(Q,P)=\int \varphi (\frac{\dr Q}{\dr P})\ dP $ if $ P<<Q$ and 

$K(Q,P)=+\infty $ otherwise,\\ where $\varphi : x\mapsto xln(x)-x+1$ is strictly convex.\\
Let us present some well-known properties of the Kullback-Leibler divergence.
\begin{property}\label{Kmini}
We have $K(P,Q)=0\Leftrightarrow P=Q.$
\end{property}
\begin{property}\label{K-SCI}
The divergence function $Q\mapsto K(Q,P)$ is convex, lower semi-continuous (l.s.c.) - for the topology that makes all the applications of the form $Q\mapsto\int fdQ$ continuous where $f$ is bounded and continuous - as well as l.s.c. for the topology of the uniform convergence.
\end{property}
\begin{property}[corollary  (1.29), page 19 of \cite{MR926905}]\label{ExitenceDeLEntropieDesProj}
If $T:(X,A)\to (Y,B)$ is measurable and if $K(P,Q)<\infty,$ then $K(P,Q)\geq K(PT^{-1},QT^{-1}),$
with equality being reached when $T$ is surjective for $(P,Q)$.
\end{property}
\begin{theorem}[theorem III.4 of \cite{AZE97}]\label{azeIII4}
Let $f:I\to \R$ be a convex function.  Then $f$ is a Lipschitz function in all compact intervals $[a,b]\subset int\{I\}.$ In particular, $f$ is continuous on $int\{I\}$.
\end{theorem}
\vskip 3mm
\noindent {\bf A.2.}{\it Useful lemmas}

\begin{lemma} \label{ProjBorne}
Let $f$ be a density in $\R^d$ bounded and positive.  Then, any projection density of $f$ - that we will name $f_a$, with $a\in\R^d_*$ - is also bounded and positive in $\R$.
\end{lemma}
\begin{lemma} \label{LoiCondBorne}
Let $f$ be a density in $\R^d$ bounded and positive.  Then any density $f(./\transp ax)$, for any $a\in\R^d_*$, is also bounded and positive.
\end{lemma}
\begin{lemma}\label{GkBorne}
If $f$ and $g$ are positive and bounded densities, then $g^{(k)}$ is positive and bounded.
\end{lemma}
\begin{lemma}\label{aleph}
Let $f$ be an absolutely continuous density, then, for all sequences $(a_n)$ tending  to $a$ in $\R^d_*$, the sequence $f_{a_n}$ uniformly converges towards $f_a$.
\end{lemma}
$\mathcal Proof$ :\\
For all $a$ in $\R^d_*$, let $F_a$ be the cumulative distribution function of $\transp aX$ and $\psi_a$ be a complex function defined by $\psi_a(u,v)=F_a(\cR e(u+iv))+iF_a(\cR e(v+iu))$, for all $u$ and $v$ in $\R$.\\First, the function $\psi_a(u,v)$ is an analytic function, because $x\mapsto f_a(\transp ax)$ is continuous and as a result of the corollary of Dini's second theorem - according to which 
{\it "A sequence of cumulative distribution functions which pointwise converges on $\R$ towards a continuous cumulative distribution function $F$ on $\R$, uniformly converges towards $F$ on $\R$"}-
we deduct that, for all sequences $(a_n)$ converging towards $a$, $\psi_{a_n}$ uniformly converges towards $\psi_a$.
Finally, the Weierstrass theorem, (see proposal $(10.1)$ page 220 of \cite{JD1980}), implies that all sequences $\psi_{a,n}'$ uniformly converge towards $\psi_a'$, for all $a_n$ tending to $a$. We can therefore conclude.\hfill$\Box$
\begin{lemma} \label{compacité-1}
The set $\Gamma_c$ is closed in $L^1$ for the topology of the uniform convergence.
\end{lemma}
\begin{lemma} \label{compacité}
For all $c>0$, we have $\Gamma_c\subset \overline B_{L^1}(f,c),$ where $B_{L^1}(f,c)=\{p\in L^1;\|f-p\|_1\leq c\}$.
\end{lemma}
\begin{lemma} \label{compacité+1}
 $G$ is closed in $L^1$ for the topology of the uniform convergence.
\end{lemma}
\begin{lemma}\label{18}
Let $H$ be an integrable function and let $C=\int\ H\ d\PP$ and $C_n=\int\ H\ d\Pn_n$,

$\ \ \ \ \ \ \ \ \ \ \ \ $then, $C_n-C=O_{\PP}(\frac{1}{\sqrt n}).$
\end{lemma}
\noindent {\bf B. Study of the sample}\label{truncSample}

Let $X_1$, $X_2$,..,$X_m$ be a sequence of independent random vectors with the same density  $f$. Let $Y_1$, $Y_2$,..,$Y_m$ be a sequence of independent random vectors with the same density  $g$.
Then, the kernel estimators $f_m$, and $f_{a,m}$ of $f$ and $f_{a}$, for all $a\in\R^d_*$, almost surely and uniformly converge since we assume that the bandwidth $h_m$ of these estimators meets the following conditions (see \cite{MR0345296}): \\
$(\mathcal Hyp)$: $h_m\searrow_m0$, $mh_m\nearrow_m\infty$, $mh_m/L(h_m^{-1})\to_m\infty$ and $L(h_m^{-1})/LLm\to_m\infty$, with $L(u)=ln(u\vee e)$.\\
Let us consider 
$A_0(m,a)=\frac{1}{m}\Sigma_{i=1}^mln\{\frac{g_{a}(\transp aY_i) }{f_{a,m}(\transp aY_i)}\}\frac{g_{a}(\transp aY_i) }{g(Y_i)},$
$A_0'(m,a)=\frac{1}{m}\Sigma_{i=1}^m(\frac{g_a(\transp aX_i) }{f_{a,m}(\transp aX_i)}-1) \frac{f_{a,m}(\transp aX_i)}{f_m(X_i)},$\\
$B_0(m,a)=\frac{1}{m}\Sigma_{i=1}^mln\{\frac{f_{a,m}(\transp aY_i)}{g_a(\transp aY_i) }\frac{g(Y_i)}{f_m(Y_i)}\}\frac{f_{a,m}(\transp aY_i)}{g_a(\transp aY_i)}$,
$B_0'(m,a)=\frac{1}{m}\Sigma_{i=1}^m(1-\{\frac{f_{a,m}(\transp aX_i)}{g_a(\transp aX_i) }\frac{g(X_i)}{f_m(X_i)}\})$.\\
Our goal is to estimate the maximum of $K(g_a,f_a)$ and the minimum of $K(g\frac{f_a}{g_a},f)$).\\
To achieve this, it is necessary for us to truncate $X_1$, $X_2$,..,$X_m$ and $Y_1$, $Y_2$,..,$Y_m$:\\
Let us consider now a sequence $\theta_m$ such that $\theta_m\to 0,$ and $y_m/\theta_m^2\to 0,$ where $y_m$ is defined through lemma \ref{KernRate} with $y_m=O_\PP(m^{-\frac{2}{4+d}})$. We will generate $f_m$ and $f_{b,m}$ from the starting sample and we select the $X_i$ and the $Y_i$ vectors such that $f_m(X_i)\geq \theta_m$  and $g(Y_i)\geq \theta_m$, for all $i$ and for all $b\in \R^d_*$ - for Huber's algorithm - and such that $f_m(X_i)\geq \theta_m$ and $g_{b}(\transp bY_i)\geq \theta_m$, for all $i$ and for all $b\in \R^d_*$ - for our algorithm. 
The vectors meeting these conditions will be called $X_1,X_2,...,X_n$ and $Y_1,Y_2,...,Y_n$.\\
Consequently, the next proposition provides us with the condition required to obtain our estimates
\begin{proposition}\label{QuotientDonneLoi}
Using the notations introduced in \cite{MR2054155} and in sections 3.1.1. and 4.1., it holds
\begin{eqnarray}
&&\sup_{a\in\R^d_*}|(A_0(n,a)-A_0'(n,a))-K(g_a,f_a)|\to0\text{ a.s.,}\label{lim1}\\
&&\sup_{a\in\R^d_*}|(B_0(n,a)-B_0'(n,a))-K(g\frac{f_a}{g_a},f)|\to0\text{ a.s.} \label{lim2}
\end{eqnarray}
\end{proposition}
\begin{remark}\label{Scott}
We can  take for $\theta_m$ the expression $m^{-\nu}$, with $0<\nu<\frac{1}{4+d}$. Moreover, 
to estimate $a_k$, $k\geq2$, we use the same procedure than the one we followed in order to find $a_1$ with $g_n^{(k-1)}$ instead of $g$ - since $g^{(k-1)}$ is unknown in this case.
\end{remark}
\vskip 3mm
\noindent{\bf C. Case study : $f$ is known}\label{fKnown}

In this Appendix, we study the case when $f$ and $g$ are known.
\vskip 3mm
\noindent{\bf C.1.}{\it  Convergence study at the  $k^{\text{th}}$ step of the algorithm:}

In this paragraph, when $k$ is less than or equal to $d$, we  show that the sequence $(\check a_k)_n$ converges towards $a_k$ and that the sequence $(\check g^{(k)})_n$ converges towards $g^{(k)}$. \\
Both $\check \gamma_n$ and $\check c_n(a)$ are M-estimators and estimate $a_k$ - see \cite{MR2054155}.
We state
\begin{proposition}\label{pConv2}
Assuming $(H'1)$ to $(H'3)$ hold. Both $\sup_{a\in\Theta}\|\check c_n(a)-a_k\|$ and $\check \gamma_n $ tends to $a_k$ a.s.
\end{proposition}
\noindent Finally, the following theorem shows us that $\check  g^{(k)}$ converges uniformly almost everywhere towards $g^{(k)}$, for any $k=1..d$.
\begin{theorem}\label{KRessultatPricipal}
Assumimg $(H'1)$ to $(H'3)$ hold. Then, $\check  g^{(k)}\to_n   g^{(k)}$ a.s. and uniformly a.e.
\end{theorem}
\noindent{\bf C.2.}{\it  Asymptotic Inference at the  $k^{\text{th}}$ step of the algorithm}

The following theorem shows that $\check g^{(k)}$ converges at the rate $O_{\PP}(n^{-1/2})$ in three differents cases, namely for any given $x$, with the  $L^1$ distance and with the Kullback-Leibler divergence:
\begin{theorem}\label{Superdiffconv}Assuming  $(H'0)$ to $(H'3)$ hold, for any $k=1,...,d$ and any $x\in\R^d$, we have
\begin{eqnarray}
&&|\check g^{(k)}(x)-g^{(k)}(x)|= O_{\PP}(n^{-1/2}),\label{diffconv1}\\
&&\int |\check g^{(k)}(x)-g^{(k)}(x)|dx=O_{\PP}(n^{-1/2}),\label{diffconv2}\\
&&|K(\check g^{(k)},f)-K(g^{(k)},f)|=O_{\PP}(n^{-1/2}).\label{diffconv3}
\end{eqnarray}
\end{theorem}
The following theorem shows that the laws of our estimators of $a_k$, namely $\check c_n(a_k)$ and $\check \gamma_n$, converge towards a linear combination of Gaussian variables.
\begin{theorem}\label{ KloiestimateursMin}
Assuming that conditions $(H'1)$ to $(H'6)$ hold, then\\
$\sqrt n\cA.(\check c_n(a_k)-a_k)\cvL \cB.\cN_d(0,\PP\|\frac{\dr}{\dr b}M(a_k,a_k)\|^2)+\cC.\cN_d(0,\PP\|\frac{\dr}{\dr a}M(a_k,a_k)\|^2)$ and\\
$\sqrt n\cA.(\check \gamma_n-a_k)\cvL \cC.\cN_d(0,\PP\|\frac{\dr}{\dr b}M(a_k,a_k)\|^2)+\cC.\cN_d(0,\PP\|\frac{\dr}{\dr a}M(a_k,a_k)\|^2)$\\
where $\cA=(\PP\frac{\dr^2}{\dr b\dr b}M(a_k,a_k)(\PP\frac{\dr^2}{\dr a_i\dr a_j}M(a_k,a_k)+\PP\frac{\dr^2}{\dr a_i\dr b_j}M(a_k,a_k)))$, \\
$\cC=\PP\frac{\dr^2}{\dr b\dr b}M(a_k,a_k)$ and
$\cB=\PP\frac{\dr^2}{\dr b\dr b}M(a_k,a_k)+\PP\frac{\dr^2}{\dr a_i\dr a_j}M(a_k,a_k)+\PP\frac{\dr^2}{\dr a_i\dr b_j}M(a_k,a_k).$
\end{theorem}
\noindent{\bf C.3.}{\it A stopping rule for the procedure}

We now assume that the algorithm does not stop after $d$ iterations. We then remark that, it still holds - for any $i>d$: 

$\bullet$  $g^{(i)}(x)=g(x)\Pi_{k=1}^i\frac{f_{ a_k}(\transp { a_k}x)}{[g_n^{(k-1)}]_{ a_k}(\transp { a_k}x)}$, with $g^{(0)}=g$.

$\bullet$  $K(g^{(0)},f)\geq K(g^{(1)},f)\geq K(g^{(2)},f)...\geq 0$.

$\bullet$ Theorems \ref{KRessultatPricipal}, \ref{Superdiffconv} and \ref{ KloiestimateursMin}.\\
Moreover, through remark \ref{criteria-H} page \pageref{criteria-H} and as explained in section 14 of \cite{MR790553}, the sequence $( K(g^{(k-1)}\frac{f_{a_k}}{g^{(k-1)}_{a_k}},f))_{k\geq 1}$ converges towards zero. Then, in this paragraph, we  show that $g^{(i)}$ converges towards $f$ in $i$. Finally, we  provide a stopping rule for this identification procedure.
\vskip 3mm
\noindent{\bf C.3.1.}{\it  Representation of $f$}

Under $(H'0)$, the following proposition shows us that the probability measure with density $g^{(k)}$ converges towards the probability measure with density $f$ : 
\begin{proposition}\label{cvl}
We have  $\lim_k g^{(k)}= f$ a.s.
\end{proposition}
\noindent{\bf C.3.2.}{\it  Testing of the criteria}

Through a test of the criteria, namely $a\mapsto K(g^{(k-1)}\frac{f_a}{g^{(k-1)}_a},f)$, we  build a stopping rule for this procedure. First, the next theorem enables us to derive the law of the  criteria.
\begin{theorem} \label{LOIDUCRITERE} 
Assuming that $(H'1)$ to $(H'3)$, $(H'6)$ and $(H'8)$ hold. Then,

$\sqrt n(Var_{\PP}(M(\check c_n(\check \gamma_n),\check \gamma_n)))^{-1/2}(\Pn_nM(\check c_n(\check \gamma_n),\check \gamma_n)-\Pn_nM(a_k,a_k)) \cvL \cN(0,I)$,\\
where $k$ represents the $k^{th}$ step of the  algorithm and with $I$ being the identity matrix in $\R^d$.
\end{theorem}
\noindent Note that $k$ is fixed in theorem \ref{LOIDUCRITERE} since $\check \gamma_n   =\ arg\inf_{a\in\Theta }\ \sup_{c\in\Theta }\ \Pn_nM(c,a)$ where $M$ is a known function of $k$ - see section 3.1.1. Thus, in the case where $K(g^{(k-1)}\frac{f_{a_k}}{g^{(k-1)}_{a_k}},f)= 0$, we obtain
\begin{corollary} \label{LOIDUCRITERE2} 
Assuming that $(H'1)$ to $(H'3)$, $(H'6)$, $(H'7)$ and $(H'8)$ hold. Then, 

$\sqrt n(Var_{\PP}(M(\check c_n(\check \gamma_n),\check \gamma_n)))^{-1/2}(\Pn_nM(\check c_n(\check \gamma_n),\check \gamma_n)) \cvL \cN(0,I).$
\end{corollary}
\noindent Hence, we propose the test of the null hypothesis 
$(H_0)$ : $K(g^{(k-1)}\frac{f_{a_k}}{g^{(k-1)}_{a_k}},f)= 0$ versus $(H_1)$ : $K(g^{(k-1)}\frac{f_{a_k}}{g^{(k-1)}_{a_k}},f)\not= 0$.
Based on this result, we stop the algorithm, then, defining $a_k$ as the last vector generated, we derive from corollary \ref{LOIDUCRITERE2} a $\alpha$-level confidence ellipsoid around $a_k$, namely

$\cE_k=\{b\in\R^d;\ \sqrt n(Var_{\PP}(M(b,b)))^{-1/2}\Pn_nM(b,b)\leq q_{\alpha}^{\cN(0,1)} \}$,\\
where $q_{\alpha}^{\cN(0,1)}$ is the quantile of a $\alpha$-level reduced centered normal distribution.\\ Consequently, the following corollary provides us with a confidence region for the above test:
\begin{corollary}\label{LOIDUCRITERE2coro}
$\cE_k$ is a confidence region for the test of the null hypothesis $(H_0)$ versus $(H_1)$.
\end{corollary}
\noindent{\bf D. Hypotheses' discussion}\label{DiscussHyp}\\
\noindent{\bf D.1.}{\it  Discussion on $(H'2)$.}

We verify this hypothesis in the case where :\\
$\bullet$ $a_1$ is the unique element of $\R^d_*$ such that $f(./\transp {a_1}x)=g(./\transp {a_1}x)$, i.e. $K(g(./\transp {a_1}x)f_{a_1}(\transp {a_1}x),f)=0$,(1)\\
$\bullet$ $f$ and $g$ are bounded and positive, (2)\\
$\bullet$ there exists a neighborhood $V$ of $a_k$ such that, for all $b$ in $V$ and for all positive real $A$, there exists $\cS>0$ such that $g(./ b'x)\leq \cS.f(./ b'x)$ with $\|x\|>A$ (3).\\
We remark that we obtain the same proof with $f$, $g^{(k-1)}$ and $a_k$.\\
First, (1) implies that $g\frac{f_{a_1}}{g_{a_1}}=f$. Hence, $0>\int ln(\frac{g}{f}\frac{f_{c}}{g_{c}})g\frac{f_{a_1}}{g_{a_1}}dx=-K(g\frac{f_{c}}{g_{c}},f)>-K(g,f)$ as a result of the very construction of $g\frac{f_{c}}{g_{c}}$.
Besides, (2) and (3) imply that there exists a neighborhood $V$ of $a_k$ such that, for all $c$ in $V$, there exists $\cS>0$ such that, for all $x$ in $\R^d$, $g(./ c'x)\leq \cS.f(./ c'x)$.\\
Consequently, we get $|M(c,a_1,x)|\leq |-K(g,f)|+|-(\frac{g(./ c'x)}{f(./c'x)}-1)|\leq K(g,f)+\cS+1$.\\
Finally, we infer the existence a neighborhood $V$ of $a_k$ such that, for all $c$ in $V$, 

$|M(c,a_k,x)|\leq H(x)=K(g,f)+\cS+1$ $(\PP -a.s.)$ with $\PP H<\infty$.
\vskip 3mm
\noindent{\bf D.2.}{\it  Discussion on $(H'3)$.}

We verify this hypothesis in the case where $a_1$ is the unique element of $\R^d_*$ such that $f(./\transp {a_1}x)=g(./\transp {a_1}x)$, i.e. $K(g(./\transp {a_1}x)f_{a_1}(\transp {a_1}x),f)=0$ - we obtain the same proof with $f$, $g^{(k-1)}$ and $a_k$.\\
\noindent{\it Preliminary $(A)$:
Shows that
$A=\{(c,x)\in\R^d_*\backslash \{a_1\}\times R^d;\ \frac{f_{a_1}(\transp {a_1}x)}{g_{a_1}(\transp {a_1}x)}>\frac{f_{c}(\transp cx)}{g_{c}(\transp cx)}\ and\ g(x)\frac{f_{c}(\transp cx)}{g_{c}(\transp cx)}> f(x)\}=\emptyset$
through a reductio ad absurdum, i.e. if we assume $A\not=\emptyset$.}\\
Thus, we have 
$f(x)=f(./\transp {a_1}x)f_{a_1}(\transp {a_1}x)=g(./\transp {a_1}x)f_{a_1}(\transp {a_1}x)>$ 
$g(./\transp {c}x)f_{c}(\transp {c}x)> f$, since $\frac{f_{a_1}(\transp {a_1}x)}{g_{a_1}(\transp {a_1}x)}\geq\frac{f_{c}(\transp cx)}{g_{c}(\transp cx)}$ implies $g(./\transp {a_1}x)f_{a_1}(\transp {a_1}x)=g(x)\frac{f_{a_1}(\transp {a_1}x)}{g_{a_1}(\transp {a_1}x)}\geq g(x)\frac{f_{c}(\transp cx)}{g_{c}(\transp cx)}=g(./\transp {c}x)f_{c}(\transp {c}x)$, i.e. $f>f$. We can therefore conclude.
$\\${\it Preliminary $(B)$:
Shows that $B=\{(c,x)\in\R^d_*\backslash \{a_1\}\times R^d;\ \frac{f_{a_1}(\transp {a_1}x)}{g_{a_1}(\transp {a_1}x)}<\frac{f_{c}(\transp cx)}{g_{c}(\transp cx)}\ and\ g(x)\frac{f_{c}(\transp cx)}{g_{c}(\transp cx)}< f(x)\}=\emptyset$
through a reductio ad absurdum, i.e. if we assume $B\not=\emptyset$.}\\
Thus, we have 
$f(x)=f(./\transp {a_1}x)f_{a_1}(\transp {a_1}x)=g(./\transp {a_1}x)f_{a_1}(\transp {a_1}x)< g(./\transp {c}x)f_{c}(\transp {c}x)< f$.\\
We can thus conclude as above.
$\\$Let us now prove  $(H'3)$.
We have $P M(c,a_1)- P M(c,a)=\int ln(\frac{g(x)f_{c}(\transp cx)}{g_{c}(\transp cx)f(x)})\{\frac{f_{a_1}(\transp {a_1}x)}{g_{a_1}(\transp {a_1}x)}-\frac{f_{c}(\transp cx)}{g_{c}(\transp cx)}\}g(x)dx.$
Moreover, the logarithm $ln$ is negative on 
$\{x\in\R^d_*;\ \frac{g(x)f_{c}(\transp cx)}{g_{c}(\transp cx)f(x)}<1\}$ and is positive on $\{x\in\R^d_*;\ \frac{g(x)f_{c}(\transp cx)}{g_{c}(\transp cx)f(x)}\geq1\}$.
Thus, the preliminary studies $(A)$ and $(B)$ show that $ln(\frac{g(x)f_{c}(\transp cx)}{g_{c}(\transp cx)f(x)})$ and $\{\frac{f_{a_1}(\transp {a_1}x)}{g_{a_1}(\transp {a_1}x)}-\frac{f_{c}(\transp cx)}{g_{c}(\transp cx)}\}$  always present a negative product. We can thus conclude, since $(c,a)\mapsto P M(c,a_1)- P M(c,a)$ is not null for all $c$ and for all $a\not=a_1$.\hfill$\Box$
\vskip 3mm
\noindent{\bf E. Proofs}
\begin{remark}\label{GkBor}
1/ $(H'0)$ - according to which $f$ and $g$ are assumed to be positive and bounded - through lemma \ref{GkBorne} (see page \pageref{GkBorne}) implies that $\check g^{(k)}_n$ and $\hat g^{(k)}_n$ are positive and bounded.\\
2/ Remark \ref{implyEstimBounded} implies that $f_n$, $g_n$, $\check g^{(k)}_n$ and $\hat g^{(k)}_n$ are positive and bounded since we consider a Gaussian kernel.
\end{remark}
\noindent{\bf Proof of propositions \ref{lemmeHuber0prop} and \ref{lemmeHuberprop}.}
Let us first study proposition \ref{lemmeHuberprop}.\\
Without loss of generality, we  prove this proposition with $x_1$ in lieu of $\transp aX$.\\
We define $g^*=gr$. We remark that $g$ and $g^*$ present the same density conditionally to $x_1$.
Indeed, $g^*_1(x_1)= \int g^*(x)dx_2...dx_d= \int r(x_1)g(x)dx_2...dx_d= r(x_1)\int g(x)dx_2...dx_d= r(x_1)g_1(x_1)$.\\
Thus, we can prove this proposition.
We have $g(.|x_1)=\frac{g(x_1,..., x_n)}{g_1(x_1)}$ and $g_1(x_1)r(x_1)$ is the marginal density of $g^*.$
Hence, $g^*$ is a density since $g^*$ is positive and since\\ $\int g^*dx=\int g_1(x_1)r(x_1)g(.|x_1)dx=\int g_1(x_1)\frac{f_1(x_1)}{g_1(x_1)}(\int g(.|x_1)dx_2..dx_d)dx_1=\int f_1(x_1)dx_1=1$. Moreover,
\begin{eqnarray}
  K(f,g^*)&=& \int f\{ln(f)-ln(g^*)\}dx,\label{dernière-1}\\
        &=& \int f\{ln(f(.|x_1))-ln(g^*(.|x_1))+ln(f_1(x_1))-ln(g_1(x_1)r(x_1))\}dx,\nonumber\\
        &=& \int f\{ln(f(.|x_1))-ln(g(.|x_1)) +ln(f_1(x_1))-ln(g_1(x_1)r(x_1))\}dx,\label{dernière}
\end{eqnarray}
as $g^*(.|x_1)=g(.|x_1)$. Since the minimum of this last equation (\ref{dernière}) is reached through the minimization of $\int f\{ln(f_1(x_1))-ln(g_1(x_1)r(x_1))\}dx=  K(f_1,g_1r)$, then property  \ref{Kmini} necessarily implies that $f_1=g_1r$, hence $r=f_1/g_1$.
Finally, we have $K(f,g)-  K(f,g^*)= \int f\{ln(f_1(x_1))-ln(g_1(x_1))\}dx=K(f_1,g_1),$
which completes the demonstration of proposition \ref{lemmeHuberprop}.\\
Similarly, if we replace $f^*=fr^{-1}$ with $f$ and $g$ with $g^*$, we obtain the proof of proposition \ref{lemmeHuber0prop}.\hfill$\Box$\\
\noindent{\bf Proof of propositions \ref{lemmeHuberModifprop} and \ref{liendeminimaxi}.}
The proof of proposition \ref{lemmeHuberModifprop} (resp. \ref{liendeminimaxi}) is very similar to the one for proposition \ref{lemmeHuberprop}, save for the fact we now base our reasoning at row \ref{dernière-1} on $K(g^*,f)= \int g^*\{ln(f)-ln(g^*)\}dx$ (resp. $\int g\{ln(g^*)-ln(f)\}dx$) instead of $K(f,g^*)= \int f\{ln(f)-ln(g^*)\}dx$.\hfill$\Box$\\
\noindent{\bf Proof of lemma \ref{ChangBasis}.}
\begin{lemma} \label{ChangBasis}
If the family $(a_i)_{i=1...d}$ is a basis of $\R^d$ then 

$g(./\transp{a_{1}}x,...,\transp{a_{j}}x)=n(\transp{a_{j+1}}x,...,\transp{a_{d}}x)=f(./\transp{a_{1}}x,...,\transp{a_{j}}x)$.
\end{lemma}
\noindent Putting $A=(a_1,..,a_d)$, let us determine $f$ in the $A$ basis.
Let us first study the function defined by $\psi:\R^d\to\R^d$, $x\mapsto(\transp{a_1}x,..,\transp{a_d}x).$
We can immediately say that $\psi$ is continuous and since $A$ is a basis, its bijectivity is obvious.
Moreover, let us study its Jacobian.
By definition, it is $J_\psi(x_1,\dotsc,x_d)= |(\frac{\partial\psi_i}{\partial x_j})_{1\leq i,j\leq d}|=|(a_{i,j})_{1\leq i,j\leq d}|=|A|\not=0$ since $A$ is a basis. We can therefore infer for any $x$ in $\R^d,$ there exists a unique $y$ $in$ $\R^d$ such that $f(x)=|A|^{-1}\Psi(y),$
i.e. $\Psi$ (resp. $y$) is the expression of $f$ (resp of $x$) in  basis $A$, namely
$\Psi(y)=\tilde n(y_{j+1},...,y_{d})\tilde h(y_{1},...,y_{j})$, with $\tilde n$ and $\tilde h$ being the expressions of $n$ and $h$ in the $A$ basis. 
Consequently, our results in the case where the family $\{a_j\}_{1\leq j\leq d}$ is the canonical basis of $\R^d$, still hold for $\Psi$ in the $A$ basis - see section 2.1.2. And then, if $\tilde g$ is the expression of $g$ in the $A$ basis, we have
$\tilde g(./y_1,...,y_{j})= \tilde n(y_{j+1},...,y_{d})=\Psi(./y_{1},...,y_{j})$, i.e. $g(./\transp{a_{1}}x,...,\transp{a_{j}}x)=n(\transp{a_{j+1}}x,...,\transp{a_{d}}x)$
$=f(./\transp{a_{1}}x,...,\transp{a_{j}}x)$.\hfill$\Box$\\
\noindent{\bf Proof of lemma \ref{toattaint}.}
\begin{lemma}\label{toattaint}
$\inf_{a\in\R^d_*}  K(g\frac{f_a}{g_a},f) \text{ is reached.}$
\end{lemma}
Indeed, let $G$ be $\{g\frac{f_a}{g_a};\ a\in\R^d_*\}$ and $\Gamma_c$ be $\Gamma_c=\{p;\   K(p,f)\leq c\}$ for all $c>0$. From lemmas \ref{compacité-1}, \ref{compacité} and \ref{compacité+1} (see page \pageref{compacité}), we get $\Gamma_c\cap G$ is a compact for the topology of the uniform convergence, if $\Gamma_c\cap G$ is not empty.
Hence, and since  property \ref{K-SCI} (see page \pageref{K-SCI}) implies that  $Q\mapsto   K(Q,P)$ is lower semi-continuous in $L^1$ for the topology of the uniform convergence, then the infimum is reached in $L^1$.
(Taking for example  $c=K(g,f),$ $\Omega$  is necessarily not empty because we always have $K(g\frac{f_a}{g_a},f)\leq K(g,f)$).\hfill$\Box$\\
\noindent{\bf Proof of lemma \ref{KernRate}.}
\begin{lemma}\label{KernRate}
For any continuous density $f$, we have
$y_m=|f_m(x)-f(x)|=O_\PP(m^{-\frac{2}{4+d}})$.
\end{lemma}
Defining $b_m(x)$ as  $b_m(x)=|E(f_m(x))-f(x)|$, we have $y_m\leq |f_m(x)-E(f_m(x))|+b_m(x)$. Moreover, 
from page 150 of \cite{MR1191168}, we derive that $b_m(x)=O_\PP(\Sigma_{j=1}^dh_j^2)$ where $h_j=O_\PP(m^{-\frac{1}{4+d}})$. Then, we infer $b_m(x)=O_\PP(m^{-\frac{2}{4+d}})$. Finally, since the central limit theorem rate  is $O_\PP(m^{-\frac{1}{2}})$, we then obtain that $y_m\leq O_\PP(m^{-\frac{1}{2}})+O_\PP(m^{-\frac{2}{4+d}})=O_\PP(m^{-\frac{2}{4+d}})$.\hfill$\Box$\\
\noindent{\bf Proof of proposition \ref{QuotientDonneLoi}.}
We  prove this proposition for $k\geq2$, i.e. in the case where $g^{(k-1)}$ is not known. The initial case using the known density $g^{(0)}=g$, will be an immediate consequence from the above. Moreover, going forward, to be more legible, we will use $g$ (resp. $g_n$) in lieu of $g^{(k-1)}$ (resp. $g^{(k-1)}_n$). We can therefore remark that we have $f(X_i)\geq \theta_n-y_n$, $g(Y_i)\geq \theta_n-y_n$ and $g_b(\transp bY_i)\geq \theta_n-y_n$, for all $i$ and for all $b\in \R^d_*$, thanks to the uniform convergence of the kernel estimators.
Indeed, we have $f(X_i)=f(X_i)-f_n(X_i)+f_n(X_i)\geq-y_n+f_n(X_i)$, by definition of $y_n$, and then  $f(X_i)\geq-y_n+\theta_n$, by hypothesis on $f_n(X_i)$. This is also true for $g_n$ and $g_{b,n}$. This entails $\sup_{b\in\R^d_*}|\frac{1}{n}\Sigma_{i=1}^{n}\{\frac{g_{b,n}(\transp bX_i)}{f_{b,n}(\transp bX_i)}-1\}\frac{f_{b,n}(\transp bX_i)}{{f_n(X_i)}}-\int\{\frac{g_b(\transp bx)}{f_b(\transp bx)}-1\}f_b(\transp bx)dx|\to 0\ a.s.$\\
Indeed, we remark that $|\frac{1}{n}\Sigma_{i=1}^{n}\{\frac{g_{b,n}(\transp bX_i)}{f_{b,n}(\transp bX_i)}-1\}\frac{f_{b,n}(\transp bX_i)}{{f_n(X_i)}}-\int\{\frac{g_b(\transp bx)}{f_b(\transp bx)}-1\}f_b(\transp bx)dx|$\\
$=|\frac{1}{n}\Sigma_{i=1}^{n}\{\frac{g_{b,n}(\transp bX_i)}{f_{b,n}(\transp bX_i)}-1\}\frac{f_{b,n}(\transp bX_i)}{{f_n(X_i)}}-\frac{1}{n}\Sigma_{i=1}^{n}\frac{g_b(\transp bX_i)}{f_b(\transp bX_i)}-1\}\frac{f_b(\transp bX_i)}{f(X_i)}$

$+\frac{1}{n}\Sigma_{i=1}^{n}\frac{g_b(\transp bX_i)}{f_b(\transp bX_i)}-1\}\frac{f_b(\transp bX_i)}{f(X_i)}-\int\{\frac{g_b(\transp bx)}{f_b(\transp bx)}-1\}f_b(\transp bx)dx|$\\
$\leq|\frac{1}{n}\Sigma_{i=1}^{n}\{\frac{g_{b,n}(\transp bX_i)}{f_{b,n}(\transp bX_i)}-1\}\frac{f_{b,n}(\transp bX_i)}{{f_n(X_i)}}-\frac{1}{n}\Sigma_{i=1}^{n}\frac{g_b(\transp bX_i)}{f_b(\transp bX_i)}-1\}\frac{f_b(\transp bX_i)}{f(X_i)}|$

$+|\frac{1}{n}\Sigma_{i=1}^{n}\frac{g_b(\transp bX_i)}{f_b(\transp bX_i)}-1\}\frac{f_b(\transp bX_i)}{f(X_i)}-\int\{\frac{g_b(\transp bx)}{f_b(\transp bx)}-1\}f_b(\transp bx)dx|$\\
Moreover, since $\int|\{\frac{g_b(\transp bx)}{f_b(\transp bx)}-1\}f_b(\transp bx)|dx\leq 2$, the law of large numbers enables us to derive:\\
$|\frac{1}{n}\Sigma_{i=1}^{n}\frac{g_b(\transp bX_i)}{f_b(\transp bX_i)}-1\}\frac{f_b(\transp bX_i)}{f(X_i)}-\int\{\frac{g_b(\transp bx)}{f_b(\transp bx)}-1\}f_b(\transp bx)dx|\to 0\ a.s..$\\
Moreover, $|\frac{1}{n}\Sigma_{i=1}^{n}\{\frac{g_{b,n}(\transp bX_i)}{f_{b,n}(\transp bX_i)}-1\}\frac{f_{b,n}(\transp bX_i)}{{f_n(X_i)}}-\frac{1}{n}\Sigma_{i=1}^{n}\frac{g_b(\transp bX_i)}{f_b(\transp bX_i)}-1\}\frac{f_b(\transp bX_i)}{f(X_i)}|$

$\ \ \ \ \ \ $ $\leq \frac{1}{n}\Sigma_{i=1}^{n}|\{\frac{g_{b,n}(\transp bX_i)}{f_{b,n}(\transp bX_i)}-1\}\frac{f_{b,n}(\transp bX_i)}{{f_n(X_i)}}-\{\frac{g_b(\transp bX_i)}{f_b(\transp bX_i)}-1\}\frac{f_b(\transp bX_i)}{f(X_i)}|$\\
and 
$|\{\frac{g_{b,n}(\transp bX_i)}{f_{b,n}(\transp bX_i)}-1\}\frac{f_{b,n}(\transp bX_i)}{{f_n(X_i)}}-\{\frac{g_b(\transp bX_i)}{f_b(\transp bX_i)}-1\}\frac{f_b(\transp bX_i)}{f(X_i)}|$
$=|\frac{g_{b,n}(\transp bX_i)-f_{b,n}(\transp bX_i)} {f_n(X_i)}-\frac{g_b(\transp bX_i)-f_b(\transp bX_i)}{f(X_i)}|$

$\leq\frac{1}{|f(X_i)|.|f_n(X_i)|}\{|f(X_i)|.|g_{b,n}(\transp bX_i)-g_b(\transp bX_i)|+|f(X_i)-f_n(X_i)|.|g_b(\transp bX_i)|$

$\ \ \ \ \ \ \ \ \ \ \ \ \ \ \ \ \ \ $ $+|f(X_i)|.|f_{b,n}(\transp bX_i)-f_b(\transp bX_i)|+|f(X_i)-f_n(X_i)|.|f_b(\transp bX_i)|\},$ 

through the introduction of terms $g_bf-g_bf$ and $ff_b-ff_b$,

$\leq\frac{O_{\PP}(1).y_n}{\theta_n.(\theta_n-y_n)}=O_{\PP}(1)\frac{1}{\frac{\theta_n^2}{y_n}-\theta_n},$
as a result of the very definitions of $\theta_n$ and $y_n$ respectively,

$\to 0,\ a.s.\text{ because, $\frac{y_n}{\theta_n^2}\to 0\ a.s.$, by hypothesis on $\theta_n$.}$\\
Consequently, $\frac{1}{n}\Sigma_{i=1}^{n}|\{\frac{g_{b,n}(\transp bX_i)}{f_{b,n}(\transp bX_i)}-1\}\frac{f_{b,n}(\transp bX_i)}{{f_n(X_i)}}-\{\frac{g_b(\transp bX_i)}{f_b(\transp bX_i)}-1\}\frac{f_b(\transp bX_i)}{f(X_i)}|\to 0$, as it is a Ces\`aro mean. This enables us to conclude. Similarly, we prove limits \ref{lim1} and \ref{lim2} page \pageref{lim1}.\hfill$\Box$\\
\noindent{\bf Proof of lemma \ref{TrucBidule}.}
\begin {lemma}\label{TrucBidule}
For any $p\leq d$, we have $f^{(p-1)}_{a_p}=f_{a_p}$ - see Huber's analytic method -, $g^{(p-1)}_{a_p}=g_{a_p}$ -  see Huber's synthetic method - and $g^{(p-1)}_{a_p}=g_{a_p}$ - see our algorithm.
\end{lemma}
$\mathcal Proof$ :\\
As it is equivalent to prove either our algorithm or Huber's, we will only develop here the proof for our algorithm.
Assuming, without any loss of generality, that the $a_i$, $i=1,..,p$, are the vectors of the canonical basis, since 
$g^{(p-1)}(x)=g(x)\frac{f_1(x_1)}{g_1(x_1)}\frac{f_2(x_2)}{g_2(x_2)}...\frac{f_{p-1}(x_{p-1})}{g_{p-1}(x_{p-1})}$ we derive immediately that $g^{(p-1)}_{p}=g_{p}$. We remark that it is sufficient to operate a change in basis on the $a_i$ to obtain the general case.
\hfill$\Box$\\
\noindent{\bf Proof of lemma \ref{imFree}.}
\begin {lemma}\label{imFree}
If there exits $p$, $p\leq d$, such that $K(g^{(p)},f)=0$, then the family  of $(a_i)_{i=1,..,p}$ - derived from the construction of $g^{(p)}$ - is free and orthogonal.
\end{lemma}
$\mathcal Proof$ :\\
Without any loss of generality, let us assume that $p=2$ and that the $a_i$ are the vectors of the canonical basis. Using a reductio ad absurdum with the hypotheses  $a_1=(1,0,...,0)$ and that $a_2=(\alpha,0,...,0)$, where  $\alpha\in\R$, we get $g^{(1)}(x)=g(x_2,..,x_d/x_1)f_1(x_1)$ and $f=g^{(2)}(x)=g(x_2,..,x_d/x_1)f_1(x_1)\frac{f_{\alpha a_1}(\alpha x_1)}{[g^{(1)}]_{\alpha a_1}(\alpha x_1)}$. Hence $f(x_2,..,x_d/x_1)=g(x_2,..,x_d/x_1)\frac
{f_{\alpha a_1}(\alpha x_1)}
{[g^{(1)}]_{\alpha a_1}(\alpha x_1)}.$\\
It consequently implies that $f_{\alpha a_1}(\alpha x_1)=[g^{(1)}]_{\alpha a_1}(\alpha x_1)$ since\\
 $1=\int f(x_2,..,x_d/x_1)dx_2...dx_d=\int g(x_2,..,x_d/x_1)dx_2...dx_d\frac
{f_{\alpha a_1}(\alpha x_1)}
{[g^{(1)}]_{\alpha a_1}(\alpha x_1)}=\frac
{f_{\alpha a_1}(\alpha x_1)}
{[g^{(1)}]_{\alpha a_1}(\alpha x_1)}$.\\
Therefore, $g^{(2)}=g^{(1)}$, i.e. $p=1$ which leads to a contradiction. Hence, the family is free.\\
Moreover, using a reductio ad absurdum we get the orthogonality. Indeed, we have \\
$\int f(x)dx=1\not=+\infty=\int n(\transp{a_{j+1}}x,...,\transp{a_{d}}x)h(\transp{a_{1}}x,...,\transp{a_{j}}x)dx$.
\hfill$\Box$\\
\noindent{\bf Proof of lemma \ref{H3}.}
\begin{lemma}\label{H3}
We have $\Theta=\{b\in\Theta\ |\ \ \int(\frac{g(x)}{f(x)}\frac{f_b(\transp bx)}{g_b(\transp bx)}-1)f(x)dx<\infty\}$.
\end{lemma}
We get the result since 
$\int\ (\frac{g(x)f_b(\transp bx)}{f(x)g_b(\transp bx)}-1)f(x)\ dx=\int\ (\frac{g(x)f_b(\transp bx)}{g_b(\transp bx)}-f(x))\ dx=0$.\hfill$\Box$\\
\noindent{\bf Proof of propositions \ref{pConv2}.}
In the same manner as in Proposition 3.4 of \cite{MR2054155}, we prove this proposition through lemma \ref{H3}.\hfill$\Box$\\
\noindent{\bf Proof of propositions \ref{KernelpConv2} and \ref{pConv1}.}
Proposition \ref{KernelpConv2} comes immediately from proposition \ref{QuotientDonneLoi} page \pageref{QuotientDonneLoi} and lemma \ref{pConv2} page \pageref{pConv2}. Similarly, we prove proposition \ref{pConv1} since both $\sup_{a\in\Theta}\|\hat b_n(a)-a_k\| $  and  $\hat \beta_n $ converge toward $a_k$ a.s. in the case where $f$ is known - see also in Appendix C, where we carry out our algorithm in the case where $f$ is known.\hfill$\Box$\\
\noindent{\bf Proof of theorem \ref{KRessultatPricipal}.}
Using lemma \ref{aleph} page \pageref{aleph} and since, for any $k$, $g^{(k)}=g^{(k-1)}\frac{f_{a_k}}{g^{(k-1)}_{a_k}}$, we prove this theorem by induction.\hfill$\Box$\\
\noindent{\bf Proof of theorems \ref{KernelKRessultatPricipal} and \ref{KernelKRessultatPricipal-H}.}
We prove the theorem \ref{KernelKRessultatPricipal} by induction.
First, by the very definition of the kernel estimator $\check g^{(0)}_n=g_n$ converges towards $g$. Moreover, the continuity of $a\mapsto f_{a,n}$ and $a\mapsto g_{a,n}$ and proposition \ref{KernelpConv2} imply that $\check g^{(1)}_n=\check g^{(0)}_n\frac{f_{a,n}}{\check g^{(0)}_{a,n}}$ converges towards $g^{(1)}$. 
Finally, since, for any $k$, $\check g^{(k)}_n=\check g^{(k-1)}_n\frac{f_{\check a_k,n}}{\check g^{(k-1)}_{\check a_k,n}}$, we conclude similarly as for $\check g^{(1)}_n$. In a similar manner, we prove theorem \ref{KernelKRessultatPricipal-H}.\hfill$\Box$\\
\noindent{\bf Proof of theorem \ref{Superdiffconv}.}\\
\noindent{\bf relationship (\ref{diffconv1}).}
We consider $\Psi_j=\{\frac{f_{\check{a_j}}(\transp{\check{a_j}}x)}{[\check g^{(j-1)}]_{\check{a_j}}(\transp{\check{a_j}}x)}-\frac{f_{a_j}(\transp{a_j}x)}{[g^{(j-1)}]_{a_j}(\transp{a_j}x)}\}$.
Since $f$ and $g$ are bounded, it is easy to prove  that from a certain rank, we get, for any given $x$ in $\R^d$

$|\Psi_j| \leq max(\frac{1}{[\check g^{(j-1)}]_{\check{a_j}}(\transp{\check{a_j}}x)},\frac{1}{[ g^{(j-1)}]_{a_j}(\transp{a_j}x)})|f_{\check{a_j}}(\transp{\check{a_j}}x)-f_{a_j}(\transp{a_j}x)|$.
\begin{remark}
First, based on what we stated earlier, for any given $x$ and from a certain rank, 
there is a constant $R>0$ independent from $n$, such that
$max(\frac{1}{[\check g^{(j-1)}]_{\check{a_j}}(\transp{\check{a_j}}x)},\frac{1}{[ g^{(j-1)}]_{a_j}(\transp{a_j}x)})\leq R=R(x)=O(1).$\\
Second, since  $\check a_k$ is an $M-$estimator of $a_k$, its convergence rate is $O_{\PP}(n^{-1/2})$.
\end{remark}
Thus using simple functions, we infer an upper and lower bound for $f_{\check{a_j}}$ and for $f_{a_j}$. We therefore reach the following conclusion:
\begin{equation}\label{PsiO}
|\Psi_j|\leq O_{\PP}(n^{-1/2}).
\end{equation}
We finally obtain:

$|\Pi_{j=1}^k\frac{f_{\check{a_j}}(\transp{\check{a_j}}x)}{[\check g^{(j-1)}]_{\check{a_j}}(\transp{\check{a_j}}x)}-\Pi_{j=1}^k\frac{f_{a_j}(\transp{a_j}x)}{[g^{(j-1)}]_{a_j}(\transp{a_j}x)}|$
$=\Pi_{j=1}^k\frac{f_{a_j}(\transp{a_j}x)}{[g^{(j-1)}]_{a_j}(\transp{a_j}x)}|\Pi_{j=1}^k\frac{f_{\check{a_j}}(\transp{\check{a_j}}x)}{[\check g^{(j-1)}]_{\check{a_j}}(\transp{\check{a_j}}x)}\frac{[g^{(j-1)}]_{a_j}(\transp{a_j}x)}{f_{a_j}(\transp{a_j}x)}-1|$.\\
Based on the relationship (\ref{PsiO}), the expression   $\frac{f_{\check{a_j}}(\transp{\check{a_j}}x)}{[\check g^{(j-1)}]_{\check{a_j}}(\transp{\check{a_j}}x)}\frac{[g^{(j-1)}]_{a_j}(\transp{a_j}x)}{f_{a_j}(\transp{a_j}x)}$ tends towards 1 at a rate of $O_{\PP}(n^{-1/2})$ for all $j$.
Consequently, $\Pi_{j=1}^k\frac{f_{\check{a_j}}(\transp{\check{a_j}}x)}{[\check g^{(j-1)}]_{\check{a_j}}(\transp{\check{a_j}}x)}\frac{[g^{(j-1)}]_{a_j}(\transp{a_j}x)}{f_{a_j}(\transp{a_j}x)}$ tends towards 1 at a rate of $O_{\PP}(n^{-1/2})$. Thus from a certain rank, we get
$|\Pi_{j=1}^k\frac{f_{\check{a_j}}(\transp{\check{a_j}}x)}{[\check g^{(j-1)}]_{\check{a_j}}(\transp{\check{a_j}}x)}-\Pi_{j=1}^k\frac{f_{a_j}(\transp{a_j}x)}{[g^{(j-1)}]_{a_j}(\transp{a_j}x)}|=O_{\PP}(n^{-1/2})O_{\PP}(1).$\\ 
In conclusion, we obtain 
$|\check g^{(k)}(x)-g^{(k)}(x)|=  g(x)|\Pi_{j=1}^k\frac{f_{\check{a_j}}(\transp{\check{a_j}}x)}{[\check g^{(j-1)}]_{\check{a_j}}(\transp{\check{a_j}}x)}-\Pi_{j=1}^k\frac{f_{a_j}(\transp{a_j}x)}{[g^{(j-1)}]_{a_j}(\transp{a_j}x)}|\leq  O_{\PP}(n^{-1/2})$.\\
\noindent{\bf relationship (\ref{diffconv2}).}
The relationship \ref{diffconv1} of theorem \ref{Superdiffconv} implies that $|\frac{\check g^{(k)}(x)}{g^{(k)}(x)}-1|=O_{\PP}(n^{-1/2})$ because, for any given $x$, $g^{(k)}(x)|\frac{\check g^{(k)}(x)}{g^{(k)}(x)}-1|=|\check g^{(k)}(x)-g^{(k)}(x)|$. Consequently, there exists a smooth function $C$ of $\R^d$ in $\R^+$ such that
$\lim_{n\to\infty}n^{-1/2}C(x)=0$ and $|\frac{\check g^{(k)}(x)}{g^{(k)}(x)}-1|\leq n^{-1/2}C(x)$, for any $x$.\\
We then have
$\int |\check g^{(k)}(x)-g^{(k)}(x)|dx=\int g^{(k)}(x)|\frac{\check g^{(k)}(x)}{g^{(k)}(x)}-1|dx
\leq\int g^{(k)}(x)C(x)n^{-1/2}dx$.\\
Moreover, $\sup_{x\in\R^d}|\check g^{(k)}(x)-g^{(k)}(x)|=\sup_{x\in\R^d}g^{(k)}(x)|\frac{\check g^{(k)}(x)}{g^{(k)}(x)}-1|=\sup_{x\in\R^d}g^{(k)}(x)C(x)n^{-1/2}\to0\ a.s.$, by theorem \ref{KRessultatPricipal}.
This implies that $\sup_{x\in\R^d}g^{(k)}(x)C(x)<\infty\ a.s.$, i.e.  $\sup_{x\in\R^d}C(x)<\infty\ a.s.$ since $g^{(k)}$ has been assumed to be positive and bounded - see remark \ref{GkBor}.\\
Thus, $\int g^{(k)}(x)C(x)dx\leq \sup C.\int g^{(k)}(x)dx=\sup C<\infty$ since $g^{(k)}$ is a density, we can therefore conclude $\int |\check g^{(k)}(x)-g^{(k)}(x)|dx\leq \sup C.n^{-1/2}=O_{\PP}(n^{-1/2}).$\hfill$\Box$\\
\noindent{\bf relationship (\ref{diffconv3}).}
We have\\
$K(\check g^{(k)},f)-K( g^{(k)},f) = \int f (\varphi(\frac{\check g^{(k)}}{f})-\varphi(\frac{g^{(k)}}{f}))dx$
$ \leq \int f\ S|\frac{\check g^{(k)}}{f}-\frac{g^{(k)}}{f}|dx=S\int|\check g^{(k)}-g^{(k)}|dx$\\
with the line before last being derived from theorem  \ref{azeIII4} page \pageref{azeIII4} and where $\varphi : x\mapsto xln(x)-x+1$ is a convex function and where $S>0$. 
We get the same expression as the one found in our Proof of Relationship (\ref{diffconv2}) section, we then obtain  $K(\check g^{(k)},f)-K( g^{(k)},f) \leq O_{\PP}(n^{-1/2})$. Similarly, we get $K( g^{(k)},f)-K(\check g^{(k)},f)\leq O_{\PP}(n^{-1/2})$. We can therefore conclude.
\hfill$\Box$\\
\noindent{\bf Proof of lemma \ref{n-FunctOf-m}.}
\begin{lemma}\label{n-FunctOf-m}
We keep the notations introduced in Appendix B. It holds $n=O(m^{\frac 1 2})$.
\end{lemma}
$\mathcal Proof$ :\\
Let us first study the Huber's case.
Let $N$ be the random variable such that\\ $N=\Sigma_{j=1}^m{\bf 1}_{\{f_m(X_j)\geq \theta_m,\ g(Y_j)\geq \theta_m\}}$. Since the events $\{f_m(X_j)\geq \theta_m\}$ and $\{ g(Y_j)\geq \theta_m\}$ are independent from one another and since $\{g(Y_j)\geq \theta_m\}\subset\{ g_m(Y_j)\geq -y_m+\theta_m\}$, we can say that

$n=m.\PP(f_m(X_j)\geq \theta_m,\ g(Y_j)\geq \theta_m)\leq m.\PP(f_m(X_j)\geq \theta_m).\PP(g_m(Y_j)\geq -y_m+\theta_m)$.\\
Consequently, let us study $\PP(f_m(X_i)\geq \theta_m)$.
Let $(\xi_i)_{i=1\ldots m}$ be the sequence such that, for any $i$ and any $x$ in $\R^d$,
$\xi_i(x)=\Pi_{l=1}^d\frac{1}{(2\pi)^{1/2}h_l}e^{-\frac{1}{2}(\frac{x_l-X_{il}}{h_l})^2}-\int \Pi_{l=1}^d\frac{1}{(2\pi)^{1/2}h_l}e^{-\frac{1}{2}(\frac{x_l-X_{il}}{h_l})^2}\ f(x)dx.$
Hence, for any given $j$ and conditionally to $X_1,$ $\ldots,$ $X_{j-1}$, $X_{j+1},$ $\ldots,$ $X_m$, the variables $(\xi_i(X_j))_{i=1\ldots m}^{i\not=j}$ are i.i.d. and centered, have the same second moment, and are such that 

$|\xi_i(X_j)|\leq \Pi_{l=1}^d\frac{1}{(2\pi)^{1/2}h_l}+\Pi_{l=1}^d\frac{1}{(2\pi)^{1/2}h_l}\int  |f(x)|dx=2.(2\pi)^{-d/2}\Pi_{l=1}^dh_l^{-1}$ since $\sup_xe^{-\frac{1}{2}x^2}\leq 1$.\\
Moreover, noting that
$f_m(x)=\frac{1}{m}\Sigma_{i=1}^m\xi_i(x)+(2\pi)^{-d/2}\frac{1}{m}\Sigma_{i=1}^m\Pi_{l=1}^dh_l^{-1}\int e^{-\frac{1}{2}(\frac{x_l-X_{il}}{h_l})^2}\ f(x)dx$, \\
we have $f_m(X_j)\geq \theta_m \Leftrightarrow\frac{1}{m}\Sigma_{i=1}^m\xi_i(X_j)+(2\pi)^{-d/2}\frac{1}{m}\Sigma_{i=1}^m\Pi_{l=1}^dh_l^{-1}\int e^{-\frac{1}{2}(\frac{x_l-X_{il}}{h_l})^2}\ f(x)dx\geq \theta_m$

$\ \ \ \ \ $ $\Leftrightarrow\frac{1}{m-1}\Sigma_{\stackrel{i=1} {i\not=j}}^m\xi_i(X_j)\geq(\theta_m-(2\pi)^{-d/2}\frac{1}{m}\Sigma_{i=1}^m\Pi_{l=1}^dh_l^{-1}\int e^{-\frac{1}{2}(\frac{x_l-X_{il}}{h_l})^2}\ f(x)dx-\frac{1}{m}\xi_j(X_j))\frac{m}{m-1}$\\
with $\xi_j(X_j)=0$. Then, defining $t$ (resp. $\varepsilon$) as $t=2.(2\pi)^{-d/2}\Pi_{l=1}^dh_l^{-1}$ (resp. \\ $\varepsilon=(\theta_m-(2\pi)^{-d/2}\Pi_{l=1}^dh_l^{-1}\frac{1}{m}\Sigma_{i=1}^m\Pi_{l=1}^d\int e^{-\frac{1}{2}(\frac{x_l-X_{il}}{h_l})^2}\ f(x)dx)\frac{m}{m-1}$), the Bennet's inequality -\cite{MR0851019} page 160 - implies that 
$\PP(\frac{1}{m-1}\Sigma_{\stackrel{i=1}{i\not=j}}^m\xi_i(X_j)\geq\varepsilon/\text{$X_1,$ $\ldots,$ $X_{j-1}$, $X_{j+1},$ $\ldots,$ $X_m$})\leq2.exp(-\frac{(m-1)\varepsilon^2}{4t^2}).$\\
Finally, since the $X_i$ are i.i.d. and since $\int (\int\Pi_{l=1}^d e^{-\frac{1}{2}(\frac{x_l-y_{l}}{h_l})^2}\ f(x)dx)f(y)dy<1$, then the law of large numbers implies that 
$\frac{1}{m}\Sigma_{i=1}^m\int\Pi_{l=1}^d e^{-\frac{1}{2}(\frac{x_l-X_{il}}{h_l})^2}\ f(x)dx\to_m\int \int\Pi_{l=1}^d e^{-\frac{1}{2}(\frac{x_l-y_{l}}{h_l})^2} f(x)f(y)dxdy$ a.s.
Consequently, since $0<\nu<\frac{1}{4+d}$ - see remark \ref{Scott} - and since $e^{-x}\leq x^{-\frac{1}{2}}$ when $x>0$, we obtain, after calculation, that, from a certain rank, $exp(-\frac{(m-1)\varepsilon^2}{4t^2})=O(m^{-\frac 1 4})$, i.e., from a certain rank, $\PP(f_m(Y_j)\geq \theta_m)=O(m^{-\frac 1 4})$. Similarly, we infer $\PP(g(Y_j)\geq \theta_m)=O(m^{-\frac 1 4})$. In conclusion, we can say that 
$n=m.\PP(f_m(X_j)\geq \theta_m).\PP(g_m(Y_j)\geq \theta_m)=O(m^{\frac 1 2})$.
Similarly, we derive the same result as above for any step of our method as well as Huber's.
\hfill$\Box$\\
\noindent{\bf Proof of theorems \ref{KernelSuperdiffconv} and \ref{KernelSuperdiffconv-H}.}
First, from lemma \ref{KernRate}, we derive that, for any $x$,\\ $\sup_{a\in\R^d_*}|f_{a,n}(\transp ax)-f_a(\transp ax)|=O_{\PP}(n^{-\frac{2}{4+d}})$. Then, let us consider $\Psi_j=\frac{f_{\check{a_j},n}(\transp{\check{a_j}}x)}{\check g^{(j-1)}_{\check{a_j},n}(\transp{\check{a_j}}x)}-\frac{f_{a_j}(\transp{a_j}x)}{g^{(j-1)}_{a_j}(\transp{a_j}x)}$, we have 
$\Psi_j=\frac{1} {\check g^{(j-1)}_{\check{a_j},n}(\transp{\check{a_j}}x)g^{(j-1)}_{a_j}(\transp{a_j}x)}$
$((f_{\check{a_j},n} (\transp{\check{a_j}}x)-f_{a_j}(\transp{a_j}x))g^{(j-1)}_{a_j}(\transp{a_j}x)+f_{a_j}.(\transp{a_j}x)(g^{(j-1)}_{a_j}(\transp{a_j}x)-\check g^{(j-1)}_{\check{a_j},n}(\transp{\check{a_j}}x)))$, i.e. $|\Psi_j|=O_{\PP}(n^{-\frac{2}{4+d}})$ since $f_{a_j}(\transp{a_j}x)=O(1)$ and $g^{(j-1)}_{a_j}(\transp{a_j}x)=O(1)$. We can therefore conclude similarly as in theorem \ref{Superdiffconv} and through lemma \ref{n-FunctOf-m}. Similarly, we derive theorem \ref{KernelSuperdiffconv-H}.\hfill$\Box$\\
\noindent{\bf Proof of theorem \ref{ KloiestimateursMin}.}
First of all, we remark that hypotheses $(H'1)$ to $(H'3)$ imply that $\check \gamma_n$ and $\check c_n(a_k)$ converge towards $a_k$ in probability.
Hypothesis $(H'4)$ enables us to derive under the integrable sign after calculation,
$\PP\frac{\dr}{\dr b}M(a_k,a_k)=\PP\frac{\dr}{\dr a}M(a_k,a_k)=0,$\\
$\PP\frac{\dr^2}{\dr a_i\dr b_j}M(a_k,a_k)=\PP\frac{\dr^2}{\dr b_j\dr a_i}M(a_k,a_k)=\int \varphi"(\frac{gf_{a_k}}{fg_{a_k}}) \frac{\dr}{\dr a_i}\frac{gf_{a_k}}{fg_{a_k}}\frac{\dr}{\dr b_j}\frac{gf_{a_k}}{fg_{a_k}}\ f\ dx,$\\
$\PP\frac{\dr^2}{\dr b_i\dr b_j}M(a_k,a_k)=-\int \varphi"(\frac{gf_{a_k}}{fg_{a_k}}) \frac{\dr}{\dr b_i}\frac{gf_{a_k}}{fg_{a_k}}\frac{\dr}{\dr b_j}\frac{gf_{a_k}}{fg_{a_k}}\ f\ dx$, $\PP\frac{\dr^2}{\dr a_i\dr a_j}M(a_k,a_k)=\int \varphi'(\frac{gf_{a_k}}{fg_{a_k}})\frac{\dr^2}{\dr a_i\dr a_j}\frac{gf_{a_k}}{fg_{a_k}}\ f\ dx,$\\ and consequently 
$\PP\frac{\dr^2}{\dr b_i\dr b_j}M(a_k,a_k)=-\PP\frac{\dr^2}{\dr a_i\dr b_j}M(a_k,a_k)=-\PP\frac{\dr^2}{\dr b_j\dr a_i}M(a_k,a_k),$
which implies,\\
$\frac{\dr^2}{\dr a_i\dr a_j}K(g\frac{f_{a_k}}{g_{a_k}},f)=\PP\frac{\dr^2}{\dr a_i\dr a_j}M(a_k,a_k)-\PP\frac{\dr^2}{\dr b_i\dr b_j}M(a_k,a_k),$

$=\PP\frac{\dr^2}{\dr a_i\dr a_j}M(a_k,a_k)+\PP\frac{\dr^2}{\dr a_i\dr b_j}M(a_k,a_k)$
$=\PP\frac{\dr^2}{\dr a_i\dr a_j}M(a_k,a_k)+\PP\frac{\dr^2}{\dr b_j\dr a_i}M(a_k,a_k).$\\
The very definition of the estimators $\check \gamma_n$ and $\check c_n(a_k)$, implies that
$
\left\{
\begin{array}{rl}
\Pn_n\frac{\dr}{\dr b}M(b,a)=0\\
\Pn_n\frac{\dr}{\dr a}M(b(a),a)=0
\end{array}
\right.
$\\
ie $\left\{
\begin{array}{rl}
\Pn_n\frac{\dr}{\dr b}M(\check c_n(a_k),\check \gamma_n)=0\\
\Pn_n\frac{\dr}{\dr a}M(\check c_n(a_k),\check \gamma_n)+\Pn_n\frac{\dr}{\dr b}M(\check c_n(a_k),\check \gamma_n)\frac{\dr}{\dr a}\check c_n(a_k)=0,
\end{array}
\right.$ i.e.
$\left\{
\begin{array}{rl}
\Pn_n\frac{\dr}{\dr b}M(\check c_n(a_k),\check \gamma_n)=0\ (E0)\\
\Pn_n\frac{\dr}{\dr a}M(\check c_n(a_k),\check \gamma_n)=0\ (E1)
\end{array}
\right.$.\\
Under $(H'5)$ and $(H'6)$, and using a Taylor development of the  $(E0)$ (resp. $(E1)$) equation, we infer there exists $(\overline c_n, \overline \gamma_n)$ (resp. $(\tilde c_n, \tilde \gamma_n)$) on the interval $[(\check c_n(a_k),\check \gamma_n),(a_k,a_k)]$ such that \\
$-\Pn_n\frac{\dr}{\dr b}M(a_k,a_k)=[\transp{(\PP\frac{\dr^2}{\dr b\dr b}M(a_k,a_k))}+o_{\PP}(1),\transp{(\PP\frac{\dr^2}{\dr a\dr b}M(a_k,a_k))}+o_{\PP}(1)]a_n.$\\
(resp. $-\Pn_n\frac{\dr}{\dr a}M(a_k,a_k)=[\transp{(\PP\frac{\dr^2}{\dr b\dr a}M(a_k,a_k))}+o_{\PP}(1),\transp{(\PP\frac{\dr^2}{\dr a^2}M(a_k,a_k))}+o_{\PP}(1)]a_n$)\\
with $a_n=(\transp{(\check c_n(a_k)-a_k)},\transp{(\check \gamma_n-a_k)})$.
Thus we get

$\sqrt n a_n=\sqrt n
\left[
\begin{array}{ccc}
\PP\frac{\dr^2}{\dr b^2}M(a_k,a_k) & \PP\frac{\dr^2}{\dr a\dr b}M(a_k,a_k)\\
\PP\frac{\dr^2}{\dr b\dr a}M(a_k,a_k) & \PP\frac{\dr^2}{\dr a^2}M(a_k,a_k) \\
\end{array}
\right]^{-1}
\left[
\begin{array}{ccc}
-\Pn_n\frac{\dr}{\dr b}M(a_k,a_k)\\
-\Pn_n\frac{\dr}{\dr a}M(a_k,a_k)\\
\end{array}
\right]+o_{\PP}(1)$

$\ \ \ $ $=\sqrt n(\PP\frac{\dr^2}{\dr b\dr b}M(a_k,a_k)\frac{\dr^2}{\dr a\dr a}K(g\frac{f_{a_k}}{g_{a_k}},f))^{-1}$\\

$\ \ \ \ \ \ \ \ $ $.\left[
\begin{array}{ccc}
\PP\frac{\dr^2}{\dr b\dr b}M(a_k,a_k)+\frac{\dr^2}{\dr a\dr a}K(g\frac{f_{a_k}}{g_{a_k}},f) & \PP\frac{\dr^2}{\dr b\dr b}M(a_k,a_k)\\
\PP\frac{\dr^2}{\dr b\dr b}M(a_k,a_k) &  \PP\frac{\dr^2}{\dr b\dr b}M(a_k,a_k)\\
\end{array}
\right].\left[
\begin{array}{ccc}
-\Pn_n\frac{\dr}{\dr b}M(a_k,a_k)\\
-\Pn_n\frac{\dr}{\dr a}M(a_k,a_k)\\
\end{array}
\right]+o_{\PP}(1)$\\
Moreover, the central limit theorem implies: 
$\Pn_n\frac{\dr}{\dr b}M(a_k,a_k)\cvL \cN_d(0,\PP\|\frac{\dr}{\dr b}M(a_k,a_k)\|^2)$,\\
$\Pn_n\frac{\dr}{\dr a}M(a_k,a_k)\cvL \cN_d(0,\PP\|\frac{\dr}{\dr a}M(a_k,a_k)\|^2)$,
since $\PP\frac{\dr}{\dr b}M(a_k,a_k)=\PP\frac{\dr}{\dr a}M(a_k,a_k)=0$, which leads us to the result.
Finally, if $f$ is known, we similarly prove theorem \ref{ KloiestimateursMax}.\hfill$\Box$\\
\noindent{\bf Proof of theorems \ref{ KernelKloiestimateursMin} and \ref{ KloiestimateursMax}.}
We get the theorem through proposition \ref{QuotientDonneLoi} and theorem \ref{ KloiestimateursMin}.\hfill$\Box$\\
\noindent{\bf Proof of proposition \ref{cvl}.}
We consider $\psi$, $\psi_a$, $\psi^{(k)}$, $\psi^{(k)}_a$ the characteristic functions of densities $f$, $f_a$, $g^{(k-1)}$ and $[g^{(k-1)}]_a$. We have
$|\psi(ta)-\psi^{(k-1)}(ta)|=|\psi_a(t)-\psi^{(k-1)}_a(t)|\leq \int |f_a(\transp ax)-[g^{(k-1)}]_a(\transp ax)|dx,$
and then
$\sup_a|\psi_a(t)-\psi^{(k-1)}_a(t)|\leq \sup_a \int |f_a(\transp ax)-[g^{(k-1)}]_a(\transp ax)|dx$\\$\leq \sup_a K([g^{(k-1)}]_a,f_a)$
since $\psi(ta)=\E(e^{it\transp ax})=\psi_a(t)$ - where $t\in\R$ and $a\in\R^d_*$ - and since the Kullback-Leibler divergence is greater than the $L^1$ distance. 
Therefore, since, as explained in section 14 of Huber's article, we have $\lim_kK([g^{(k-1)}]_{a_k},f_{a_k})=0$ we then get $\lim_k g^{(k)}=f$ - which is the Huber's representation of $f$.
Moreover, we have
$|\psi(t)-\psi^{(k)}(t)|\leq \int |f(x)-g^{(k)}(x)|dx\leq K(g^{(k)},f).$ As explained in section 14 of Huber's article and through remark \ref{criteria-H} page \pageref{criteria-H} as well as through the additive relationship of proposition \ref{lemmeHuber0prop}, we infer that $\lim_kK(g^{(k-1)}\frac{f_{a_k}}{[g^{(k-1)}]_{a_k}},f)=0$.
Consequently, we get $\lim_k g^{(k)}=f$ - which is our representation of $f$.\\
\noindent{\bf Proof of lemmas \ref{FromSection14Huber} and \ref{KernelHuber0}.}
We apply our algorithm between $f$ and $g$. There exists a sequence of densities $(g^{(k)})_k$ such that
$0=K(g^{(\infty)},f)\leq ..\leq K(g^{(k)},f)\leq..\leq K(g,f)$, (*)\\
where $g^{(\infty)}=\lim_kg^{(k)}$ which is a density by construction.
Moreover, let $(g^{(k)}_n)_k$ be the sequence of densities such that $g^{(k)}_n$ is the kernel estimate of $g^{(k)}$. Since we derive from remark \ref{GkBor} page \pageref{GkBor} an integrable upper bound of $g^{(k)}_n$, for all $k$, which is greater than $f$ - see also the definition of $\varphi$ in the proof of theorem \ref{limnk} -, then the dominated convergence theorem implies that, for any $k$, $\lim_nK(g^{(k)}_n,f_n)=K(g^{(k)},f)$, i.e., from a certain given rank $n_0$, we have 

$0\leq..\leq K(g^{(\infty)}_n,f_n)\leq ..\leq K(g^{(k)}_n,f_n)\leq..\leq K(g_n,f_n)$, (**)\\
Consequently, through lemma \ref{Ineq***} page \pageref{Ineq***}, there exists a $k$ such that

$0\leq..\leq K(\Psi^{(\infty)}_{n,k},f_n)\leq..\leq K(g^{(\infty)}_n,f_n)\leq ..\leq K(\Psi^{(\infty)}_{n,k-1},f_n)\leq..\leq K(g_n,f_n)$, (***)\\
where $\Psi^{(\infty)}_{n,k}$ is a density such that $\Psi^{(\infty)}_{n,k}=\lim_kg^{(k)}_n$.
Finally, through the dominated convergence theorem and taking the limit as $n$ in (***) we get
$0=K(g^{(\infty)},f) = \lim_nK(g^{(\infty)}_n,f_n)\geq \lim_nK(\Psi^{(\infty)}_{n,k},f_n)\geq 0$.
The dominated convergence theorem enables us to conclude:

$0=\lim_nK(\Psi^{(\infty)}_{n,k},f_n)=\lim_n\lim_kK(g^{(k)}_n,f_n)$. Similarly, we get lemma \ref{KernelHuber0}.\hfill$\Box$\\
\noindent{\bf Proof of lemma \ref{Ineq***}.}
\begin{lemma}\label{Ineq***}
Keeping the notations of the proof of lemma \ref{FromSection14Huber}, we have 

$0\leq..\leq K(\Psi^{(\infty)}_{n,k},f_n)\leq..\leq K(g^{(\infty)}_n,f_n)\leq ..\leq K(\Psi^{(\infty)}_{n,k-1},f_n)\leq..\leq K(g_n,f_n)$, (***)
\end{lemma}
$\mathcal Proof$ :\\
First, as explained in section 4.2., we have $K(f^{(k)},g)-K(f^{(k+1)},g)=K(f^{(k)}_{a_{k+1}},g_{a_{k+1}})$. Moreover, through remark \ref{criteria-H} page \pageref{criteria-H}, we also derive that $K(f^{(k)},g)=K(g^{(k)},f)$. Then, $K(f^{(k)}_{a_{k+1}},g_{a_{k+1}})$ is the decreasing step of the relative entropies in (*) and leading to $0=K(g^{(\infty)},f)$. Similarly, the very construction of (**), implies that $K(f^{(k)}_{a_{k+1},n},g_{a_{k+1},n})$ is the decreasing step of the relative entropies in (**) and leading to  $K(g^{(\infty)}_n,f_n)$.
Second, through the conclusion of the section 4.2. and the lemma 14.2 of Huber's article, we obtain that $K(f^{(k)}_{a_{k+1},n},g_{a_{k+1},n})$ converges - decreasingly and in $k$ - towards a positive function of $n$ - that we will call $\xi_n$.
Third, the convergence of $(g^{(k)})_k$ - see proposition \ref{cvl} - implies that, for any given $n$, the sequence $(K(g^{(k)}_n,f_n))_k$ is not finite. Then, through relationship $(**)$, there exists a $k$ such that $0<K(g^{(k-1)}_n,f_n)-K(g^{(\infty)}_n,f_n)<\xi_n$.\\
Consequently, since $Q\mapsto K(Q,P)$ is l.s.c. - see property \ref{K-SCI} - relationship (**) implies (***).
\hfill$\Box$\\
\noindent{\bf Proof of theorems \ref{limnk} and \ref{limnk-H}.}
We recall that $g_n^{(k)}$ is the kernel estimator of $\check g^{(k)}$.
Since the Kullback-Leibler divergence is greater than the $L^1$-distance, we then have
$\lim_n\lim_k K(g_n^{(k)},f_n)\geq \lim_n\lim_k \int |g_n^{(k)}(x)-f_n(x)|dx$.
Moreover, the Fatou's lemma implies that

$\lim_k \int |g_n^{(k)}(x)-f_n(x)|dx\geq \int \lim_k \big[|g_n^{(k)}(x)-f_n(x)|\big]dx=\int |[\lim_k g_n^{(k)}(x)]-f_n(x)|dx$ and\\
$\lim_n \int |[\lim_k g_n^{(k)}(x)]-f_n(x)|dx\geq \int \lim_n \big[|[\lim_k g_n^{(k)}]-f_n|\big]dx=\int |[\lim_n\lim_k g_n^{(k)}(x)]-\lim_nf_n(x)|dx$.\\
We then obtain that $0=\lim_n\lim_k K(g_n^{(k)},f_n)\geq \int |\lim_n\lim_k g_n^{(k)}(x)-\lim_nf_n(x)|dx\geq0$, i.e. that $\int |\lim_n\lim_k g_n^{(k)}(x)-\lim_nf_n(x)|dx=0$. Moreover, for any given $k$ and any given $n$, the function $g_n^{(k)}$ is a convex combination of multivariate Gaussian distributions. As derived at remark  \ref{implyEstimBounded}, for all $k$, the determinant of the covariance of the random vector - with density $g^{(k)}$ - is greater than or equal to the product of a positive constant times the determinant of the covariance of the random vector  with density $f$.
The form of the kernel estimate therefore implies that there exists an integrable function $\varphi$ such that, for any given $k$ and any given $n$, we have $|g_n^{(k)}|\leq \varphi$. 
Finally, the dominated convergence theorem enables us to say that 
$\lim_n\lim_k g_n^{(k)}=\lim_n f_n=f$, since $f_n$ converges towards $f$ and since $\int |\lim_n\lim_k g_n^{(k)}(x)-\lim_nf_n(x)|dx=0$. Similarly, we get theorem \ref{limnk-H}.\hfill$\Box$\\
\noindent{\bf Proof of theorem \ref{LOIDUCRITERE}.}
Through a Taylor development of $\Pn_nM(\check c_n(a_k),\check \gamma_n)$ of rank 2, we get at point $(a_k,a_k)$: 
$\Pn_nM(\check c_n(a_k),\check \gamma_n)$
$=\Pn_nM(a_k,a_k)+\Pn_n\frac{\dr}{\dr a}M(a_k,a_k)\transp {(\check \gamma_n-a_k)}+\Pn_n\frac{\dr}{\dr b}M(a_k,a_k)\transp {(\check c_n(a_k)-a_k)}$

$\ \ \ \ \ \ $ $+\frac{1}{2}\{\transp {(\check \gamma_n-a_k)}\Pn_n\frac{\dr^2}{\dr a\dr a}M(a_k,a_k)(\check \gamma_n-a_k)+\transp {(\check c_n(a_k)-a_k)}\Pn_n\frac{\dr^2}{\dr b\dr a}M(a_k,a_k)(\check \gamma_n-a_k)$

$\ \ \ \ \ \ $ $+\transp {(\check \gamma_n-a_k)}\Pn_n\frac{\dr^2}{\dr a\dr b}M(a_k,a_k)(\check c_n(a_k)-a_k)+\transp {(\check c_n(a_k)-a_k)}\Pn_n\frac{\dr^2}{\dr b\dr b}M(a_k,a_k)(\check c_n(a_k)-a_k)\}$\\
Thus, lemma \ref{18} implies $\Pn_nM(\check c_n(a_k),\check \gamma_n)=\Pn_nM(a_k,a_k)+O_{\PP}(\frac{1}{n}),$
\\ i.e. $\sqrt n(\Pn_nM(\check c_n(a_k),\check \gamma_n)-\PP M(a_k,a_k))=\sqrt n(\Pn_nM(a_k,a_k)-\PP M(a_k,a_k))+o_{\PP}(1).$ \\Hence $\sqrt n(\Pn_nM(\check c_n(a_k),\check \gamma_n)-\PP M(a_k,a_k))$ abides by the same limit distribution as\\ $\sqrt n(\Pn_nM(a_k,a_k)-\PP M(a_k,a_k))$, which is $\cN(0,Var_{\PP}(M(a_k,a_k)))$.\hfill$\Box$\\
\noindent{\bf Proof of theorems \ref{KernelLOIDUCRITERE} and \ref{KernelLOIDUCRITERE-H}.}
Through proposition \ref{QuotientDonneLoi} and theorem \ref{LOIDUCRITERE}, we derive theorem \ref{KernelLOIDUCRITERE}. Similarly, we get theorem \ref{KernelLOIDUCRITERE-H}.\hfill$\Box$

\end{document}